\documentclass[12pt, letter]{article}

\usepackage{graphicx}
\usepackage{amsmath, amssymb,theorem,verbatim}
\usepackage{apacite}
\usepackage{setspace}
\usepackage{array}
\usepackage{multirow}
\usepackage{booktabs}

\newcommand\VRule[1][\arrayrulewidth]{\vrule width #1}

\long\def\comment#1{}

\newcommand{\bG}{\mathbf{G}}

\newcommand{\bI}{\mathbf{I}}

\newcommand{\bZ}{\mathbf{Z}}

\newcommand{\Ex}{\mathbb{E}}
\newcommand{\cov}{\textrm{cov}}
\newcommand{\var}{\textrm{var}}
\newcommand{\cum}{\textrm{cum}}
\newcommand{\Dcon}{\stackrel{D}{\rightarrow}}

\newcommand{\Pcon}{\stackrel{\mathcal{P}}{\rightarrow}}
\newtheorem{theorem}{Theorem}[section]

\newtheorem{lemma}{Lemma}[section] 

\newtheorem{defin}{Definition}[section]

\newtheorem{assumption}{Assumption}[section]

\newtheorem{remark}{Remark}[section]

\begin{document}

\title{{\bf The quantile spectral density and comparison based tests for 
nonlinear time series}}
\author{Junbum Lee and Suhasini Subba Rao\\
Department of Statistics, Texas A\&M University \\
College Station, U.S.A.\\
{\tt jlee@stat.tamu.edu} and {\tt suhasini@stat.tamu.edu}} 
\date{}

\maketitle

\begin{abstract}
In this paper we consider tests for nonlinear time series, which are
motivated by the notion of serial dependence. 
The proposed tests are based on
comparisons with the quantile spectral density, which can be
considered as a quantile version of the usual spectral density
function. The quantile spectral density `measures' the sequential dependence
structure of a time series, and is well defined under relatively weak
mixing conditions. We propose an estimator for the quantile spectral
density and derive its asympototic sampling properties. 
We use the quantile spectral density to construct a
goodness of fit test for time series and explain how this test can
also be used for comparing the sequential dependence structure of two
time series. The asymptotic sampling properties of the test statistic is
derived under the null and an alternative. Furthermore, 
a bootstrap procedure it proposed to obtain a finite sample approximation. 
The method is illustrated with simulations and some real data examples.


\end{abstract}

\section{Introduction}

The analysis of most time series is based on a set of assumptions,
which in practice need to be tested. This is usually done through a
goodness of fit test. 
The majority of goodness of fit tests for time series are based on fitting the 
conjectured model to the data, estimating the residuals of the
model and testing for lack of correlation, normally with a Ljung-Box
type test  (see for example, \citeA{p:and-93} and \citeA{p:hon-96}).  
\citeA{p:che-deo-04} propose a test based on the spectral density,  
and \citeA{p:hal-pur-92} propose robust tests 
based on ranks.  If one restricts the class of models to just 
linear time series models, then these tests can correctly identify 
the model. However, problems can arise, if one widens the class of
models and allow for nonlinear time series. For example, if the time
series were to satisfy an ARCH process, then it will be uncorrelated,
but it is not independent. Furthermore, the squares will satisfy an
autoregressive representation, with errors which are martingale
differences. Therefore, correlation based test for nonlinear time series models
may not identify the model. 

\citeA{p:neu-pap-08} propose a goodness of fit test for Markov
time series models based on the one step ahead transition distribution. 
But this test is specifically for Markov models.
An alternative approach is to generalise the notion of correlation to
measuring the general dependence between pairs of random variables in
a time series. This notion is usually called serial dependence, 
and can be traced dates back to \citeA{p:hoe-48}.      
\citeA{p:ska-tjo-93} and \citeA{p:hon-00} use this definition to  
test for serial independence of a time series. \citeA{p:hon-98} takes
these notions further, and  generalises the
spectral density to sequential dependence. He does this by defining
the generalised spectral density, which is the Fourier transform of the
characteristic function of pair-wise dependent data. He uses this
device in \citeA{p:hon-98} and \citeA{p:hon-lee-03} to test for goodness
of fit of a time series model, mainly through the analysis of the
estimated residuals.  
%The advantage of this approach
%is that, in general, under the null hpothesis, the test statistic is asymptotically
%`distribution free', in other words auxillary parameters do not need to be estimated. 
However, sometimes the residuals cannot be or are not easy to estimate. 
For example, it is possible to estimate the residuals of an ARCH $(X_{t}=Z_{t}\sigma_{t})$, 
possible but difficult with a GARCH and usually impossible for many models of the type 
$X_{t} = g(X_{t-1},\varepsilon_{t})$. 

In this paper, we use the notion of serial dependence to test for
goodness of fit, but without estimating the residuals. 
In Section \ref{sec:motivate} we motivate our test by
considering the Microsoft daily log return data and compare it with the
GARCH$(1,1)$ model (one of the standard models for such 
data sets). We show that though the GARCH model seems to model well
some of the stylised facts of this data, ie. the
uncorrelatedness of the data and positive correlation in the absolute values, 
if one made a deeper analysis and compared the correlation
of other transformations such as $\cov(I(X_{t}\leq x),I(X_{t+r}\leq
y))$ (where $I$ denotes the indictator function), there is large difference
between the data and  GARCH model. This motivates us to define the 
{\it quantile} autocovariance function and the {\it quantile} spectral
density. The quantile spectral density can be considered as a measure
of serial dependence of a time series.  
In Sections \ref{sec:quantile} and \ref{sec:test} we propose a method for 
estimating the quantile
spectral density, and use the quantile spectral density as the basis of a
test which compares the quantile
spectral density estimator with the spectral density estimator under
the null hypothesis. The asymptotic sampling properties of the
quantile spectral density estimator are derived in Section
\ref{sec:sampling-spec}. Recently there have been several
articles defining and estimating the spectral density of sequential
dependence. In particular, 
\citeA{p:li-08}, \citeA{p:hag-11} and \citeA{p:det-11}
define spectral density functions similar to the quantile spectral
density, however these authors, estimate the periodogram and the
quantile spectral density using $L_{1}$ methods. In contrast, our
approach is motivated by the definition of the periodogram, this leads to
an estimator of the quantile spectral density with an analytic form, thus can easily be used in both goodness of 
fit and other tests. However, 
it is interesting, and rather surprising, to note that the $L_{1}$ 
estimator proposed in
\citeA{p:det-11} and our estimator of the quantile spectral density
share similar asymptotic properties.  
In Section \ref{sec:sampling} we derive the asymptotic
sampling properties of the test statistic. An advantage of our approach
is that it can easily be extended to test other quantities, for 
example with a small
adaption it can be used to test for equality of
serial dependence of two time series, this is considered in Section
\ref{sec:more}. In Section \ref{sec:practical} we propose a bootstrap
method for estimating the finite sampling distribution of the test
statistic under the null. 
The proofs can be found in the Appendix and technical report.

\section{The quantile spectral density and the test statistic}

\subsection{Motivation}\label{sec:motivate}

\begin{figure}[h!]
\begin{center}
\includegraphics[width=14cm, height=8cm]{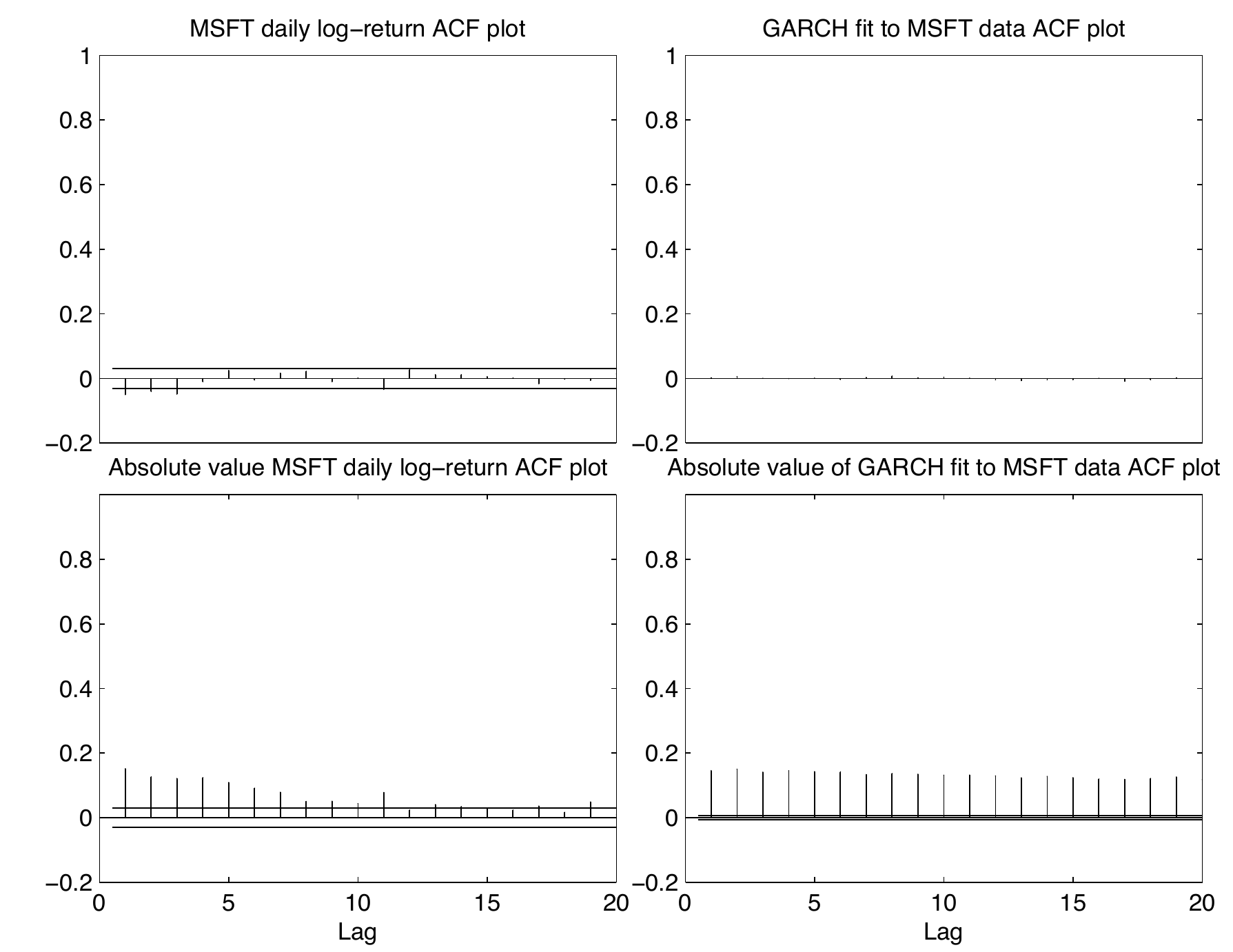}
\end{center}
\caption{\label{fig:1}The ACF plots of $\{X_{t}\}$ and $\{|X_{t}|\}$
  of the MSFT and the corresponding GARCH model}
\end{figure}

\begin{figure}[h!]
\begin{center}
\includegraphics[width=14cm,
height=6cm]{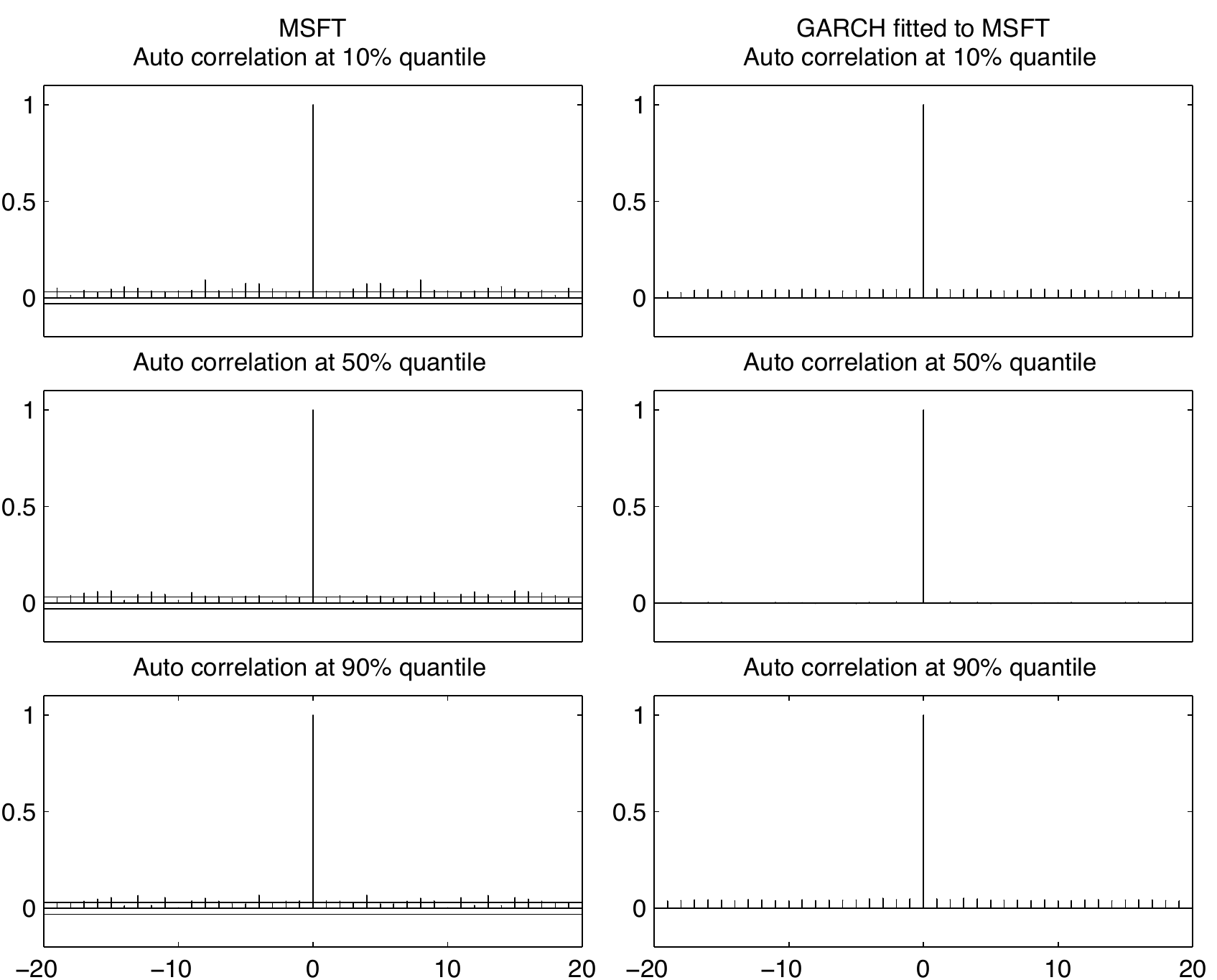}
\includegraphics[width=14cm,
height=6cm]{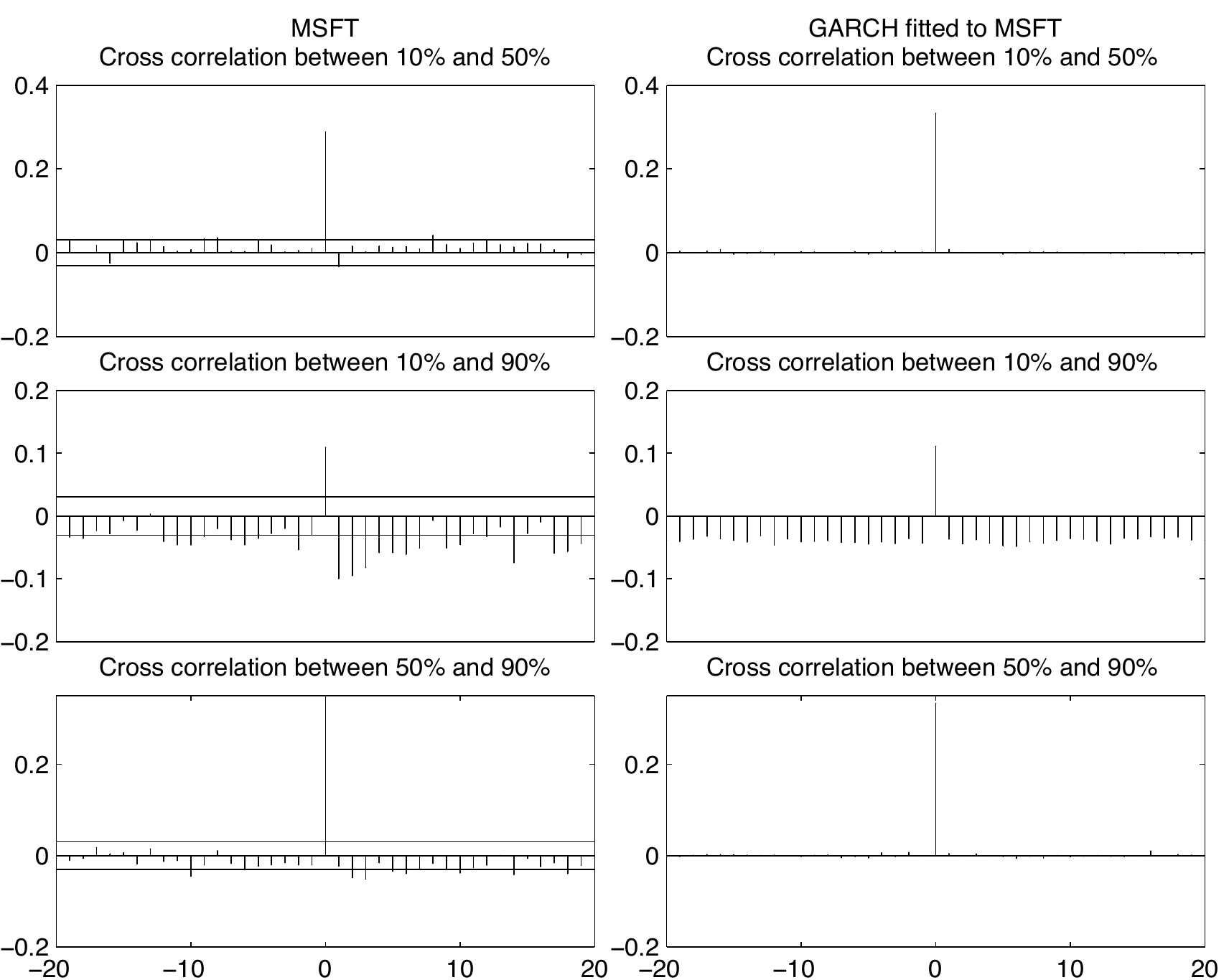}
\end{center}
\caption{\label{fig:2} The quantile covariance of the MSFT and the
  corresponding GARCH}
\end{figure}

To motivate our approach, we analyze the Microsoft daily log returns
(MSFT) between March 1986 - June 2003, 
which we denote as $\{X_{t}\}$. One argument for fitting GARCH types models to 
financial data is their
ability to model the so called `stylised facts' seen in such data sets.  
We now demonstrate why this is the case for the MSFT data
(see \citeA{b:ziv-09}). Using the 
maximim likelihood, the GARCH model which best fits the log
differences of the MSFT is 
$X_{t} = \mu + \varepsilon_t$, where 
$\epsilon_t = \sigma_{t}Z_{t}$, $\sigma_{t}^{2} = 
a_{0} + a_{1}\varepsilon_{t-1}^{2}+b\sigma_{t-1}^{2}$ ($\{Z_{t}\}$ are
independent, identically distributed standard normal random variables), with 
$\mu = 1.56\times 10^{-3}$, 
$a_{0} = 1.03\times 10^{-5}$, $a_{1} = 0.06$ and $b=0.925$.  In Figure
\ref{fig:1} we give the sample  autocorrelation plots of $\{X_{t}\}$ and
$\{|X_{t}|\}$, together with the autocorrelation plots of the corresponding 
GARCH$(1,1)$ model. 
Comparing the two plots, it appears that the GARCH$(1,1)$ captures 
the `stylised facts' in the Microsoft data, such as
the near zero autocorrelation of the observations and the persistant
positive autocorrelations of the absolute log
returns. However, if we want to check the suitability of the GARCH model 
for modelling 
the general pair-wise dependence structure, that is the joint distribution of 
$(X_{s},X_{t})$ for all $s$ and $t$ (often called {\it sequential dependence}), then 
we need to look beyond the covariance of $\{X_{t}\}$ and $\{|X_{t}|\}$. To make a 
more general comparison we transform the data into indicator variables 
$\{I(X_{t}\leq x)\}$ and check the correlation structure of the 
indicator variables
over various $x$. For example, define the multivariate vector time series
$\underline{Y}_{t} = (I(X_{t}\leq q_{0.1}), I(X_{t}<q_{0.5}) ,
I(X_{t}\leq q_{0.9}))$,  where $q_{\alpha}$ denotes the estimated 
$\alpha$-percentile of $X_{t}$.   
%(note that Stoffer and Wall (1991)\cite{}, 
%also transformed a univariate categorical time series into a
%multivariate time series, in order to define the 
%spectral envelope. 
Plots of the cross-covariances of $\underline{Y}_{t}$ and 
the corresponding GARCH model (with Gaussian innovations)
are given in Figure \ref{fig:2}. 
In Figure \ref{fig:2}, there are clear differences in the
dependence structure of the data and the GARCH model. The $10$th,
$50$th and $90$th percentiles correspond to large negative, 
zero and large positive values of $X_{t}$ (big negative change, 
no change and large positive changes in the returns). 
In order to do the analysis, we will use the following observations.
By using that $\cov(I(X_{0}\leq x),I(X_{r}\leq y)) = P(X_{0}\leq x,X_{r}\leq y) - 
P(X_{0}\leq x)P(X_{r}\leq y)$, for all  $x,y\in \mathbb{R}$ we have 
\begin{eqnarray*}
  \cov(I(X_{0}\leq x),I(X_{r}\leq y)) =  \cov(I(X_{0}\geq
  x),I(X_{r}\geq y)) = -\cov(I(X_{0}\leq x),I(X_{r} >  y)).
\end{eqnarray*}
From Figure \ref{fig:2} we observe:
\begin{itemize}
\item The ACF of $I(X_{t}\leq q_{0.5})$ of the GARCH is zero. This is
  due to the symmetry of the GARCH process: given the event $X_{0}\leq
  0$, we have equal chance $X_{r}> 0$ and $X_{r} < 0$ (ie.  
$\cov(I(X_{0}\leq 0),I(X_{r}\leq 0))= -\cov(I(X_{0}\leq 0),I(X_{r}
>0))$). This means that $\cov(I(X_{0}\leq 0),I(X_{r}\leq 0)) = 0$. On
the other hand, for the MSFT data we see that there is a clear positive 
correlation in the
sample autocorrelation of $\{I(X_{t}< 0)\}$. One interpretation for this
behaviour, is that a decrease in consecutive values, is likely to lead to
future decreases. 
\item The cross correlation $\cov(I(X_{0}< q_{0.1}),I(X_{r} < q_{0.9}))$,
where $\{X_{t}\}$ comes from a GARCH is symmetric, ie. 
$\cov(I(X_{0} < {q_{0.1}}),
I(X_{r} < q_{0.9})) = \cov(I(X_{0} < q_{0.1}),I(X_{-r} <
q_{0.9}))$. On the other hand, the corresponding sample
cross-correlations of the MSFT is not symmetric. Thus the GARCH
process is time reversible, whereas it appears that the MSFT data may
not be.    
\end{itemize}  
The cross and autocovariances in Figure \ref{fig:2} are a 
graphical representation of the serial dependence structure of the time series.  
These plots suggest that for the MSFT data the GARCH model may not be the most 
appropriate model, especially if validity is based on modelling the serial 
dependence structure.  In the sections below we will test this.  

\subsection{The quantile spectral density function}\label{sec:quantile}

We now formalise the discussion above.  Let us suppose that $\{X_{t}\}$ is a
strictly stationary time series. It is obvious that the cross
covariance of the indicator functions $\{I(X_{t} \leq x),I(X_{t}\leq
y)\}$ is  
\begin{eqnarray*}
C_{r}(x,y)&:=&\cov(I(X_{0}\leq x),I(X_{r}\leq y))= P(X_{0}\leq
x,X_{r}\leq y) - P(X_{0}\leq x)P(X_{r}\leq y).  
\end{eqnarray*}
\citeA{p:ska-tjo-93} and \citeA{p:hon-00}
use a similar quantity to test for serial independence of a time series.
We call $C_{r}(\cdot)$ the {\it quantile covariance}. If $\{X_{t}\}$ is an
$\alpha$-mixing time series with mixing rate $s >
1$ ($s$ is defined in Assumption \ref{assum:A}, below) it can be
shown that $\sup_{x,y}\sum_{r}|\cov(I(X_{0}\leq x),I(X_{r}\leq
y))|<\infty$, thus for all $x,y\in \mathbb{R}$,  it's Fourier transform 
\begin{eqnarray*}
G(x,y;\omega) &=& \frac{1}{2\pi}\sum_{r}C_{r}(x,y;\omega)\exp(ir\omega),
\end{eqnarray*}
is well defined. Since $G(x,y;\omega)$ can be considered as the
cross-spectral density of $\{I(X_{t}< x),I(X_{t}< y)\}$, 
we call $G(\cdot)$ the {\it quantile spectral density}. 

%\begin{remark}
%Under the stronger condition of $\psi$-mixing of the time series $\{X_{t}\}$ 
%(where the rate $s > 1$), 
%it can be shown that $|G(x,y;\omega)|\leq P(X_{0}\leq x)P(X_{0}\leq y)$. 
%\end{remark}

\subsubsection{Properties of the quantile spectral density}

The quantile spectral density carries all the information about the serial 
dependence structure of the time series. 
For example (i) if $\{X_{t}\}$ is serially independent, then
$G$ does not depend on $\omega$ and $G(x,y;\omega) = \rho(x,y)$ 
(ii) if for all $r$, the distribution function of $(X_{0},X_{r})$
is identical to the distribution function of $(X_{0},X_{-r})$, then 
$G(\cdot)$ will be real and (iii) for any given $x$ and $y$, $G$ gives information
about any periodicities that may exists at a given threshold. In addition,
$G(\cdot)$ captures the covariance structure of any transformation
of $\{X_{t}\}$ . For example, consider the transformation  $\{h(X_{t})\}$, 
then it is straightforward to show that the spectral density of the time series
$\{h(X_{t})\}$ is 
\begin{eqnarray*}
f_{h}(\omega) = \frac{1}{2\pi}\sum_{r}\cov(h(X_{0}),h(X_{r}))\exp(ir\omega)
 = \int\int h(x)h(y)G(dx,dy;\omega). 
\end{eqnarray*} 
Of course, $G(x,y;\omega)$ only captures the serial dependency, and 
may miss higher order structure. Only in the case that 
$\{X_{t}\}$ is Markovian, does $G(x,y;\omega)$ capture the
entire joint distribution of $\{X_{t}\}$.  

\begin{remark}
The quantile spectral density is closely related to the generalised spectral density
introduced in \citeA{p:hon-98}. He defines the 
generalised spectral density as $h(x,y;\omega) =
\sum_{r}\cov(\exp(ixX_{0}),\exp(iyX_{r}))\exp(ir\omega)$. 
Essentially, this is the Fourier transform of the  
characteristic function of pairwise distributions minus their marginals, 
therefore the relationship between the 
quantile spectral density and the generalised spectral density is
analogous to that between the distribution function and the
characteristic function of a random variable. Hong (1998, 2003)
uses the generalised spectral
density as a tool in various tests goodness of fit tests, which are mainly based 
on the residual. On the other hand, the goodness of fit test that we propose, is based on  
checking for similarity between the estimated quantile spectral density and the
proposed spectral density. 
\end{remark}

\begin{remark}[The Copula spectral density]
A closely related quantity to the quantile spectral density is the
copula spectral density, which is defined as 
\begin{eqnarray}
\label{eq:copula}
G_{\mathcal{C}}(u_{1},u_{2};\omega) &=&
\frac{1}{2\pi}\sum_{r}\mathcal{C}_{r}(u_{1},u_{2};
\omega)\exp(ir\omega), 
\end{eqnarray}
where $\mathcal{C}_{r}(u_{1},u_{2}) = \cov(I(F(X_{0}) \leq
u_{1}),I(F(X_{r})\leq u_{2})) = \Ex(I(F(X_{0}) \leq
u_{1})I(F(X_{r})\leq u_{2})) - u_{1}u_{2}$, and $F(\cdot)$ is
marginal distribution function of $\{X_{t}\}$. Note that by definition
$u_{1},u_{2}\in [0,1]$. Thus, unlike the quantile spectral density,
the copula spectral density is invariant to any monotonic
transformation of $\{X_{t}\}$, for example mean and variance shifts.
By considering the ranks of $\{X_{t}\}$, the methods detailed in the section below can also be used to estimate 
$G_{\mathcal{C}}$. 
\citeA{p:det-11} have recently proposed $L_{1}$-methods for estimating
$G_{\mathcal{C}}$, and the asymptotic sampling properties
have been derived for this estimator. 
\end{remark}

\begin{figure}[h!]
\begin{center}
\includegraphics[width=14cm, height=6cm]{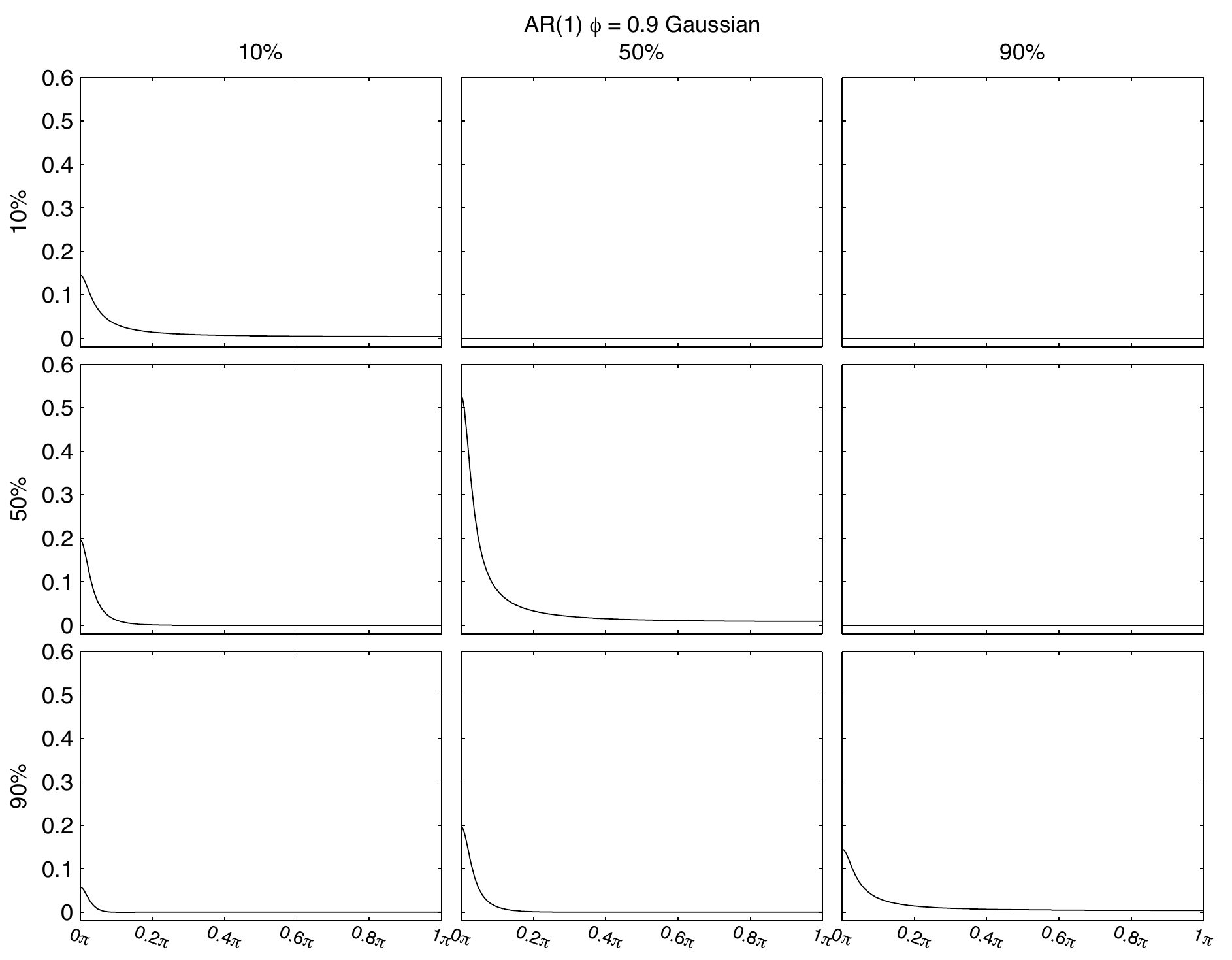}
\end{center}
\caption{\label{fig:3} The quantile spectral density of $X_{t} = 0.9X_{t-1}+Z_{t}$}
\end{figure}

\begin{figure}[h!]
\begin{center}
\includegraphics[width=14cm, height=6cm]{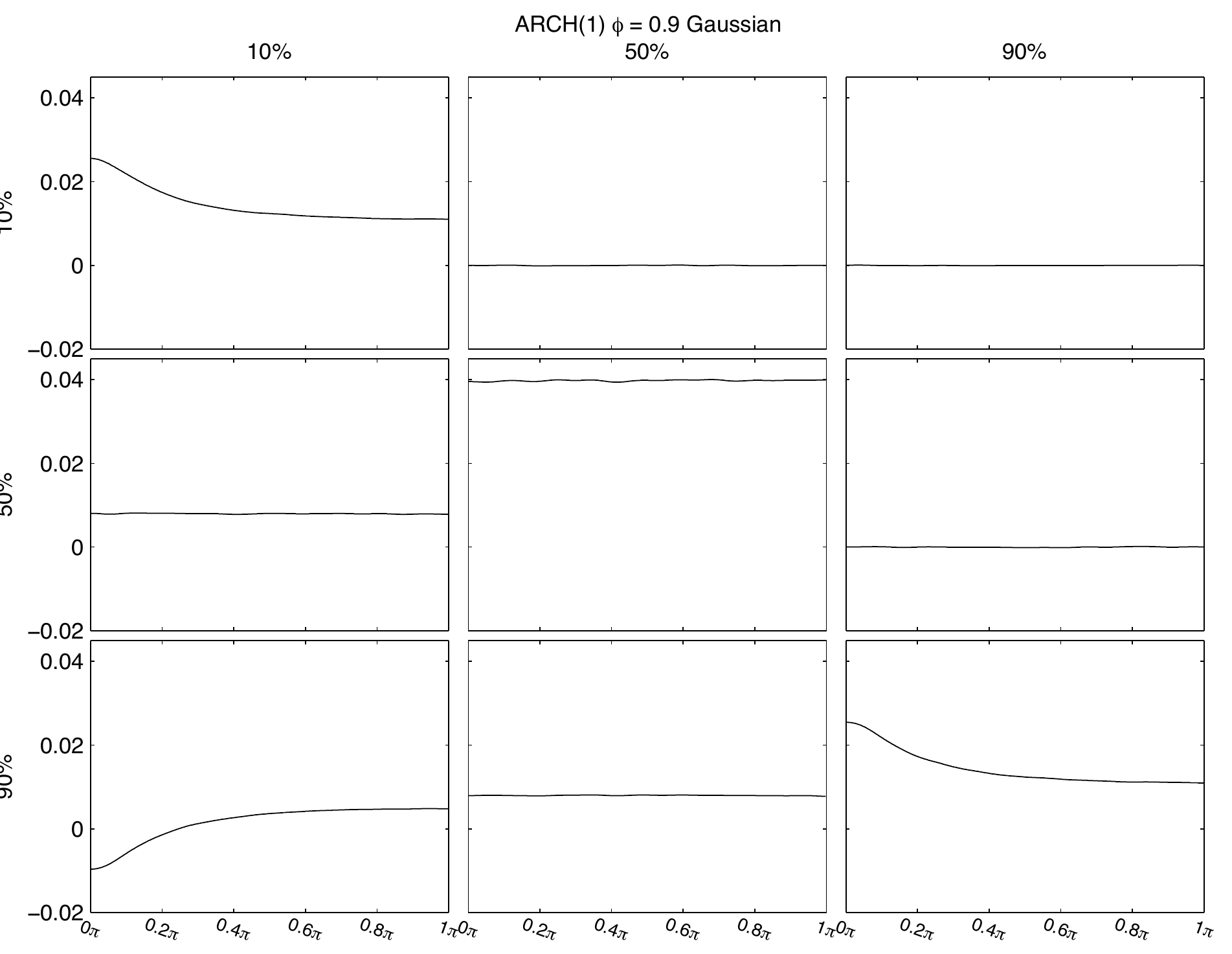}
\end{center}
\caption{\label{fig:4}The quantile spectral density of
  $X_{t}=\sigma_{t}Z_{t}$, where $\sigma_{t}^{2} = 1/1.9 + 0.9X_{t-1}^{2}$}
\end{figure}

\begin{figure}[h!]
\begin{center}
\includegraphics[width=14cm, height=6cm]{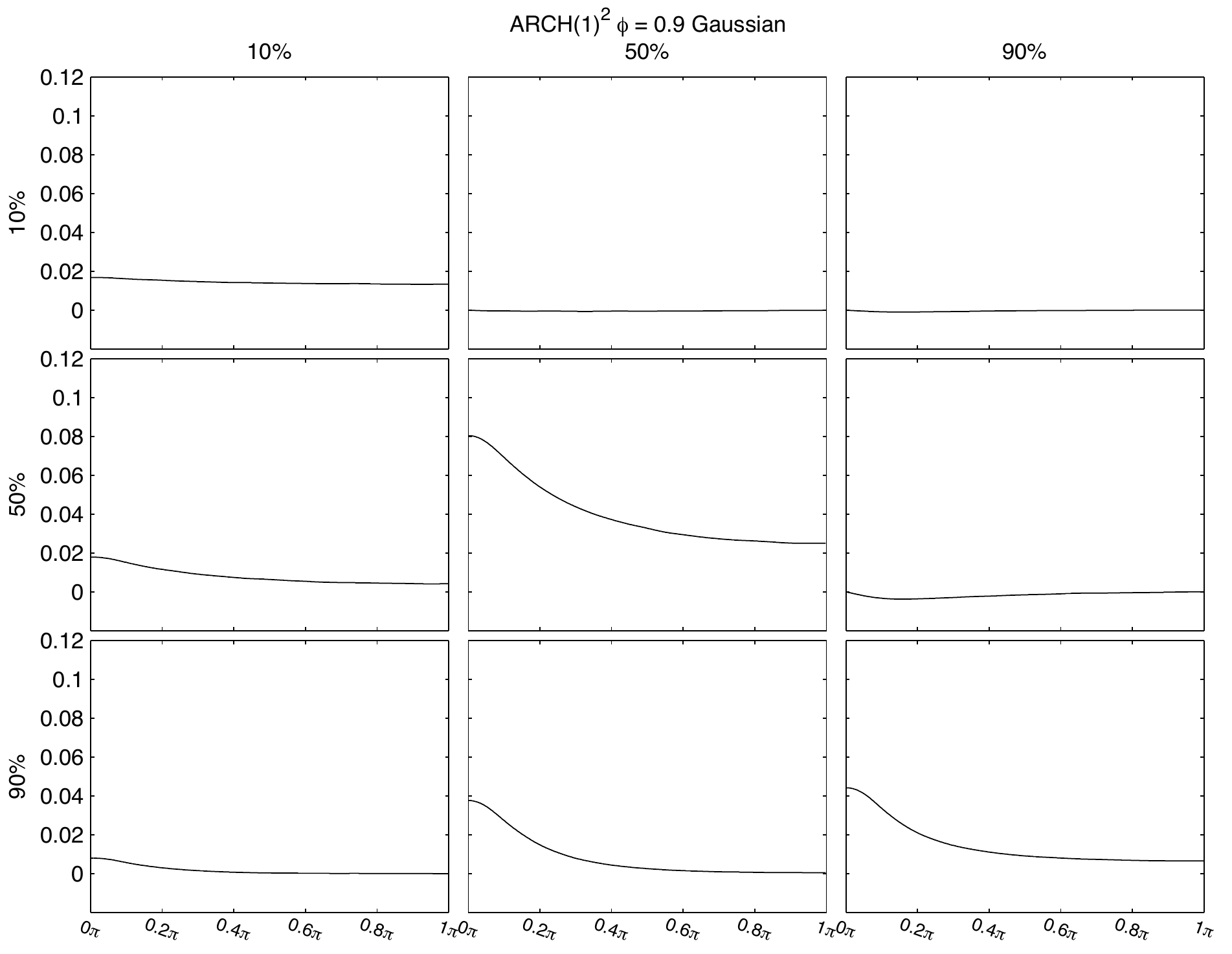}
\end{center}
\caption{\label{fig:5} The quantile spectral density of
  $X_{t}^{2}=\sigma_{t}^{2}Z_{t}^{2}$, where $\sigma_{t}^{2} = 1/1.9 + 0.9X_{t-1}^{2}$}
\end{figure}

In Figures \ref{fig:3}, \ref{fig:4} and \ref{fig:5} we plot 
the quantile spectral density for the autoregressive ($X_{t} = 0.9X_{t-1}+Z_t$),
ARCH ($X_{t}  =\sigma_{t} Z_{t}$ with 
$\sigma_{t}^{2} = 1/1.9 + 0.9 X_{t-1}^{2}$) and squared ARCH, with independent, identically distributed  (iid) Gaussian
innovations $Z_t$. The diagonals are of $G(x,x;\omega)$, the 
lower triangle contains the real part of $G(x,y;\omega)$ and the upper triangle 
the imaginary part of $G(x,y;\omega)$.  
We observe that the AR and ARCH quantile spectral densities are very
different. The AR has a similar shape for all $x$, whereas for
the ARCH, it is flat (like the spectral density of uncorrelated data)
at about the $50\%$ percentile, but moves away from flatness at the
extremes. Furthermore, recalling that the AR and ARCH squared have the same
spectral density (if the moments of the ARCH squared exists), there is
a large difference between the quantile spectral density of the AR and
the ARCH squared.

\subsubsection{Estimating the quantile spectral density}

The quantile spectral density $G(x,y;\omega)$ can be considered as the 
cross spectral density of the bivariate time series $\{I(X_{t}\leq
x),I(Y_{t}\leq y)\}$. Therefore, our estimator of $G(x,y;\omega)$ is
motivated by the classical cross spectral. To do this we define the 
class of lag windows we shall use. 
\begin{defin}\label{def:lag}
The lag window takes the form 
\begin{eqnarray*}
\lambda(u) = \big(\sum_{j=-r}^{r}a_{r}\exp(i2\pi ru) - 
\sum_{j=1}^{r}b_{j}|u|^{j}\big)I_{[-1,1]}(u), 
\end{eqnarray*}
where $I_{[-1,1]}(u) = 1$ for $u\in [-1,1]$ and zero otherwise. 
This class of lag windows is quite large, and includes the 
truncated window, the 
Bartlett window and general Tukey window (see, for example, 
\citeA{b:pri-88} Section 6.2.3 for properties of these lag windows). 
\end{defin}
To obtain an estimator of $G$, we define the centralised, transformed variable 
$Z_{t}(x) = I(X_{t}\leq x) - \widehat{F}_{T}(x)$ (where $\widehat{F}_{T}(x) =
\frac{1}{T}\sum_{t}I(X_{t}\leq x)$). 
We estimate the quantile covariance $C_{r}(x,y) = P(X_{0}\leq x,X_{r}\leq y) - 
P(X_{0}\leq x)P(X_{r}\leq y)$ with
$\widehat{C}_{r}(x,y) = \frac{1}{T}\sum_{}Z_{t}(x)Z_{t+r}(y)$, and use as
an estimator of $G$
\begin{eqnarray}
\widehat{G}_{T}(x,y;\omega_{k}) &=&
\frac{1}{2\pi}\sum_{r}\lambda_{M}(r)\widehat{C}_{r}(x,y)\exp(ir\omega_{k}) \label{eq:spectral-est}\\
 &=&  \sum_{s}K_{M}(\omega_{k} - \omega_{s})J_{T}(x;\omega_{s})
\overline{J_{T}(y;\omega_{s})}, \nonumber
\end{eqnarray}
where $\lambda_{M}(r) = \lambda(r/M)$,
$J_{T}(x;\omega) = \frac{1}{\sqrt{2\pi T}}
\sum_{t=1}^{T}Z_{t}(x)\exp(it\omega)$ and $K_{M}(\omega) = 
\frac{1}{T}\sum_{r}\lambda_{M}(r)\exp(ir\omega)$.

\subsection{The test statistic}\label{sec:test}

The proposed test is based on the fit of the estimated quantile spectral
density to the conjectured quantile spectral density. More precisely, 
we test  $H_{0}: G(x,y;\omega)=G_{0}(x,y;\omega)$ against 
$H_{A}: G(x,y;\omega) \neq G_{0}(x,y;\omega)$, where $G$ is the quantile spectral
density of $\{X_{t}\}$, $G_{0}(x,y;\omega) =
\frac{1}{2\pi}\sum_{r}C_{0,r}(x,y)\exp(ir\omega)$ 
and  $C_{0,r}(x,y) = F_{0,r}(x,y) - F_{0}(x)F_{0}(y)$. Thus under the 
null the marginal distribution is $F_{0}(\cdot)$ and the
joint distribution is $F_{0,r}(\cdot)$. We use the quadratic distance
to measure the distance between the estimated quantile
spectral density and the conjectured spectral density, and define the 
test statistic as   
\begin{eqnarray}
\mathcal{Q}_{T} &=& 
\frac{1}{T}\sum_{k=1}^{T}\int|\widehat{G}_{T}(x,y;\omega_{k}) - 
\frac{1}{2\pi}\sum_{r}\lambda_{M}(r)C_{0,r}(x,y)\exp(ir\omega_{k})|^{2}dF_{0}(x)dF_{0}(y)\nonumber\\
&=& \frac{1}{T}\sum_{k=1}^{T}\int|\widehat{G}_{T}(x,y;\omega_{k}) - 
\sum_{s=1}^{T}K_{M}(\omega_{k}-\omega_{s})G_{0}(\omega_{s})|^{2}dF_{0}(x)dF_{0}(y)\nonumber\\
 &=& \frac{1}{2\pi}\sum_{r}\lambda_{M}(r)^{2}\int\int
\big|\hat{C}_{r}(x,y) - C_{0,r}(x,y)\big|^{2}dF_{0}(x)dF_{0}(y), \label{eq:testQ}
\end{eqnarray}
where the above immediately follows from Parseval's theorem. 
The choice of lag window will have an influence on the type of
alternatives the test can detect. 
For example, the truncated window ($\lambda(u) = I_{[-1,1]}(u)$) gives equal 
weights to all the quantile covariances, whereas the Bartlett window  
($\lambda(u) = (1-|u|)I_{[-1,1]}(u)$) gives more weight to the lower order lags.  
Therefore the tests ability to detect the alternative will 
depend on which lags of 
the quantile covariance deviates the most from the null, 
and the weight the lag window 
places on these. We derive the asymptotic distribution of
$\mathcal{Q}_{T}$ in Section \ref{sec:sampling}. 

\begin{remark}
The test can be adapted to be invariant to monotonic
transformations (such as shifts of mean and variance). This can be done by replacing
the quantile spectral density with the
copula spectral density $G_{\mathcal{C}}(\cdot)$ defined
in (\ref{eq:copula}). In this case the null is $H_{0}:G_{\mathcal{C}}(x,y;\omega)
= G_{\mathcal{C},0}(x,y;\omega) = \frac{1}{2\pi}\sum_{r}\mathcal{C}_{0,r}(u_{1},u_{2};
\omega)\exp(ir\omega)$ against   $H_{A}:G_{\mathcal{C}}(x,y;\omega)
\neq G_{\mathcal{C},0}(x,y;\omega)$. The test statistic in this case
is 
\begin{eqnarray*}
\mathcal{Q}_{T,\mathcal{C}} &=& 
\frac{1}{T}\sum_{k=1}^{T}\int|\widehat{G}_{T,\mathcal{C}}(u_{1},u_{2};\omega_{k}) - 
\frac{1}{2\pi}\sum_{r}\lambda_{M}(r)\mathcal{C}_{0,r}(u_{1},u_{2})
\exp(ir\omega_{k})|^{2}du_{1}du_{2}, 
\end{eqnarray*}
where we estimate
$\widehat{G}_{T,\mathcal{C}}(u_{1},u_{2};\omega_{k})$ in the same way
as we have estimated $\widehat{G}_{T}$ in (\ref{eq:spectral-est}) but 
replace $\{X_{t}\}_{t}$ with $\{\hat{F}_{T}(X_{t})\}_{t}$.  The
distribution of $\mathcal{Q}_{T,\mathcal{C}}$ is beyond the scope of
the current paper. 
\end{remark}

\section{Sampling properties}\label{sec:sampling}

In this section we derive the sampling properties of the quantile
spectral density $\widehat{G}_{T}$ and the test statistic
$\mathcal{Q}_{T}$. 
We will use the $\alpha$-mixing assumptions below. 
\begin{assumption}\label{assum:A}
Let us suppose that $\{X_{t}\}$ is a strictly stationary $\alpha$-mixing 
time series such that
\begin{eqnarray*}
%\label{eq:assum-mix}
\sup_{\substack{A\in \sigma(X_{r},X_{r+1},\ldots) \\
B\in \sigma(X_{0},X_{-1},\ldots)}}|P(A\cap B) - P(A)P(B)| \leq \alpha(r),
\end{eqnarray*}
where $\alpha(r)$ are the mixing coefficients which satisfy 
$\alpha(r)\leq K|r|^{-s}$ for some $s > 2$.
\end{assumption}

\subsection{Sampling properties of $\widehat{G}_{T}$}\label{sec:sampling-spec}

In the following lemma we derive the limiting distribution of $\widehat{G}_{T}$, this 
will allow us to construct point wise confidence intervals for $G$. 
\begin{theorem}\label{theorem:QS}
Suppose Assumption \ref{assum:A} holds. Then 
\begin{eqnarray*}
\Ex(\widehat{G}_{T}(x,y;\omega)) = G(x,y;\omega) + O(\frac{1}{M^{s-1}}), 
\end{eqnarray*}
and for $0 < \omega_{k} < \pi$ we have 
\begin{eqnarray*}
V_{T}(x,y;\omega_{k})^{-1/2}
\left(
\begin{array}{c}
\Re \widehat{G}_{T}(x,y;\omega_{k}) - \Re\Ex(\widehat{G}_{T}(x,y;\omega_{k})) \\
\Im \widehat{G}_{T}(x,y;\omega_{k}) -  \Im\Ex(\widehat{G}_{T}(x,y;\omega_{k}))
\end{array}
\right) &\Dcon& 
\mathcal{N}\big(0,I_{2} \big) \\
V_{T}(x,x;\omega_{k})^{-1/2}
\left(\widehat{G}_{T}(x,x;\omega_{k}) - \Ex(\widehat{G}_{T}(x,x;\omega_{k}))\right) 
&\Dcon& \mathcal{N}(0,1), 
\end{eqnarray*}
where $M\rightarrow \infty$ and  $M/T\rightarrow 0$ as $T\rightarrow \infty$, 
\begin{eqnarray*}
V_{T}(x,y;\omega_{k}) = 
\sum_{k=1}^{T}K_{M}(\omega_{k} - \omega_{s})^{2}\left(
\begin{array}{cc}
A(x,y;\omega_{s}) & C(x,y;\omega_{s}) \\
C(x,y;\omega_{s}) & B(x,y;\omega_{s}) \\
\end{array}
\right) = O(\frac{M}{T}), 
\end{eqnarray*}
and
\begin{eqnarray*}
 A(x,y;\omega_{s}) &  = &  \frac{1}{2}\bigg(G(x,x;\omega_{s})G(y,y;\omega_{s}) + \Re
G(x,y;\omega_{s})^{2} - \Im G(x,y;\omega_{s})^{2}\bigg) \\
B(x,y;\omega_{s}) & = & \frac{1}{2}\bigg(G(x,x;\omega_{s})
G(y,y;\omega_{s}) + \Im G(x,y;\omega_{s})^{2} - \Re
G(x,y;\omega_{s})^{2} \bigg)\\
C(x,y;\omega_{s}) & = &\Re G(x,y;\omega_{s})\Im G(x,y;\omega_{s}).
\end{eqnarray*}
\end{theorem}
Thus, if $\frac{M}{T} >> \frac{1}{M^{2(s-1)}}$, in other words the variance of 
$\widehat{G}_{T}$ dominates the bias, then we can use the above result to construct
confidence intervals for $G$. 

\subsection{Sampling properties of test statistic under the null hypothesis}\label{sec:sampling}

We now derive the limiting distribution of the test statistic under
the null hypothesis. 
Let 
\begin{eqnarray}
\label{eq:ETstar}
E_{T} &=& \frac{1}{T}\int \int W_{M}(\omega- \theta)^{2}
G(x,x;\theta)G(y,y;\theta)dF_{0}(x)dF_{0}(y)d\theta d\omega \nonumber\\
V_{T} &=& \frac{4}{T^{2}}\int \int
\Delta_{M}(\theta_{1} - \theta_{2})^2\prod_{i=1}^{2}
G(x_{1},x_{2};\theta_{i})G(y_{1},y_{2};\theta_{i})d\theta_{i}
dF_{0}(x_{i})dF_{0}(y_{i}),
\end{eqnarray}
where 
\begin{eqnarray}
W_M(\theta) = \frac{T}{2\pi} K_M(\theta) = \frac{1}{2\pi}
\sum_{r} \lambda_M(r) \exp(ir\theta) \nonumber \\
\label{eq:delta}
\Delta_{M}(\theta_{1}-\theta_{2}) = 
\int W_{M}(\omega - \theta_{1})W_{M}(\omega - \theta_{2})d\omega.
\end{eqnarray}

\begin{lemma}\label{lemma:mean-var}
Suppose that Assumption \ref{assum:A} holds and $G(\cdot)$ is the
quantile spectral density of $\{X_{t}\}$. Then under the null
hypthesis we have 
\begin{eqnarray*}
\Ex\big(\mathcal{Q}_{T}\big) = E_{T} + O(\frac{1}{T}) = O(\frac{M}{T}) 
\textrm{ and } 
\var\big(\mathcal{Q}_{T} \big) = V_{T} + O(\frac{1}{T})  = O(\frac{M}{T^{2}}).
\end{eqnarray*}
\end{lemma}
Using the above we obtain the limiting distribution under the null. 
\begin{theorem}\label{theorem:null}
Suppose that Assumption \ref{assum:A} holds. Then under the null
hypthesis we have
\begin{eqnarray*}
V_{T}^{-1/2}\big(\mathcal{Q}_{T} - E_{T}\big)\Dcon \mathcal{N}(0,1)
\end{eqnarray*}
as $M\rightarrow \infty$ and $M/T\rightarrow 0$ as $T\rightarrow \infty$. 
\end{theorem}
Using estimates of $\widehat{G}_{T}(\cdot)$, 
$E_{T}$ and $V_{T}$ can both be estimated. Thus by using the above result, 
we reject the null at the $\alpha \%$ level if
$V_{T}^{-1/2}\big(\mathcal{Q}_{T} - E_{T}\big) > z_{1-\alpha}$ 
(where $z_{1-\alpha}$ denotes the $1-\alpha$ quantile of a standard normal 
distribution). 

\subsection{Behaviour of the test statistic under the alternative hypothesis}

We now examine the behaviour of the test statistic under the alternative
$H_{A}:G(x,y;\omega) = G_{1}(x,y;\omega) = \frac{1}{2\pi}
\sum_{r}\big(F_{r,1}(x,y) -
F_{1}(x)F_{1}(y)\big)\exp(ir\omega)$. To obtain the limiting distribution we decompose 
the test statistic $\mathcal{Q}_{T}$ as 
$\mathcal{Q}_{T}  = \mathcal{Q}_{T,1} + \mathcal{Q}_{T,2} + \mathcal{Q}_{T,3}$,  
where 
\begin{eqnarray*}
\mathcal{Q}_{T,1} &=& \frac{1}{T}\sum_{k=1}^{T}\int 
\big|\widehat{G}_{T}(x,y;\omega_{k}) -\Ex(\widehat{G}_{T}(x,y;\omega_{k}))\big|^{2}
dF_{0}(x)dF_{0}(y) \\
\mathcal{Q}_{T,2} &=& \frac{2}{T}\Re\sum_{k=1}^{T}\int 
\big[\hat{G}_{T}(x,y;\omega_{k}) - \Ex(\hat{G}_{T}(x,y;\omega_{k}))\big]
\big[\Ex(\hat{G}_{T}(x,y;\omega_{k})) -
\tilde{G}(x,y;\omega_{k})\big]dF_{0}(x)dF_{0}(y) \\
\mathcal{Q}_{T,3} &=& \frac{1}{T}\sum_{k=1}^{T}\int 
\big|\Ex(\hat{G}_{T}(x,y;\omega_{k})) -
\tilde{G}(x,y;\omega_{k})\big|^{2}dF_{0}(x)dF_{0}(y),
\end{eqnarray*}
and 
\begin{eqnarray*}
\tilde{G}(x,y;\omega_{k}) = \frac{1}{2\pi}
\sum_{r}\lambda_{M}(r) C_{r,0}(x,y)\exp(ir\omega) = \sum_{s}
K_M(\omega_{k} - \omega_{s}) G_{0}(x,y;\omega_{s}).
\end{eqnarray*}

From the decomposition of $\mathcal{Q}_{T}$, we observe that 
there are two stochastic terms
$\mathcal{Q}_{T,1}$ and $\mathcal{Q}_{T,2}$, and a deterministic term
$\mathcal{Q}_{T,3}$. By using Lemma \ref{lemma:mean-var}, it can be shown that
$\mathcal{Q}_{T,1} = O_{p}(\frac{M^{1/2}}{T} + \frac{M}{T})$. On the other hand,
we show in
the proof of the theorem below that
$\mathcal{Q}_{T,2}$, is of
lower order than $\mathcal{Q}_{T,1}$ and, thus, determines the distribution
of $\mathcal{Q}_{T}$. To understand the role that $\mathcal{Q}_{T,3}$
plays in the test, we replace $\tilde{G}(x,y;\omega)$ and 
$\Ex(\widehat{G}_{T}(x,y;\omega))$ with $G_{0}$ and $G_{1}$ respectively and obtain 
\begin{eqnarray*}
\mathcal{Q}_{T,3} &=&  \frac{1}{T}\sum_{k=1}^{T}\int 
\big|G_{1}(x,y;\omega_{k}) - 
G_{0}(x,y;\omega_{k})\big|^{2}dF_{0}(x)dF_{0}(y) + O(\frac{1}{M^{s-1}}).
\end{eqnarray*}
Thus $\mathcal{Q}_{T,3}$ measures the deviation of the alternative
from the null hypothesis, and shifts the mean of the test
statistic.

\begin{theorem}\label{theorem:alternative}
Suppose that Assumption \ref{assum:A} holds, and for all $r$, 
$\sup_{x,y}|C_{0,r}(x,y)| \leq K|r|^{-(2+\delta)}$, for some $\delta > 0$. 
%and $G(\cdot)$ denotes the
%quantile covariance of $\{X_{t}\}$. 
Under the alternative hypothesis we have 
\begin{eqnarray}
\label{eq:Q2}
\sqrt{T}\mathcal{Q}_{T,2} \Dcon \mathcal{N}(0,V_{T,2}),
\end{eqnarray} 
and 
\begin{eqnarray}
  \label{eq:Q3}
  \sqrt{T}\big(\mathcal{Q}_{T} - \mathcal{Q}_{T,3}\big)\Dcon 
\mathcal{N}(0,V_{T,2})
\end{eqnarray}
where $M\rightarrow \infty$ and $\sqrt{M}/T\rightarrow 0$ as
$T\rightarrow \infty$, 
%\begin{eqnarray*}
%E_{T,2} =\int\int 
%\big|G_{1}(x,y;\omega)) -G_{0}(x,y;\omega)\big|^{2}dF_{0}(x)dF_{0}(y)d\omega. 
%\end{eqnarray*}
and 
\begin{eqnarray*}
V_{T,2} &=& \frac{8}{T}\Re\int\int \Lambda_{T}(x_{1},y_{1};\omega)
\overline{\Lambda_{T}(x_{2},y_{2};\omega)} \\
&&\bigg\{G_{1}(x_{1},x_{2};\omega) G_{1}(y_{1},y_{2};\omega)  + 
G_{1}(x_{1},y_{2};\omega)
G_{1}(y_{1},x_{2};\omega)\bigg\}d\omega\prod_{i=1}^{2}dF_{0}(x_{i})dF_{0}(y_{i}) \\
&& + \frac{8}{T}\Re\int\int \Lambda_{T}(x_{1},y_{1};\omega_{1})
\overline{\Lambda_{T}(x_{2},y_{2};\omega_{2})} 
 G_{(x_{1},y_{1},x_{2},y_{2})}(\omega_{1},-\omega_{1},\omega_{2})\prod_{i=1}^{2}dF_{0}(x_{i})dF_{0}(y_{i})d\omega_{i},
\end{eqnarray*}
where $\Lambda_{T}(x,y;\omega_{s}) = 
\frac{1}{2\pi}\sum_{r}\lambda_{M}(r)^{2}(\frac{T-|r|}{T})[C_{1,r}(x,y)
-C_{0,r}(x,y)]\exp(ir\omega_{k})$ and $G_{(x_{1},y_{1},x_{2},y_{2})}$ is
the cross tri-spectral density of $\{(I(X_{t}\leq x_{1}),I(X_{t}\leq
y_{1}),I(X_{t}\leq x_{2}),I(X_{t}\leq y_{2}))\}_{t}$. 
\end{theorem}

The theorem above tells us that the mean of the test statistic is
shifted the further the alternative is from the null. Interestingly, 
we observe from the definition of $\Lambda_{T}(\cdot)$, that 
the variance also depends on the difference between the null and
alternative. However, for a fixed alternative, the
power of the test converges to $100\%$ as the sample size grow. 

\section{Testing for equality of serial dependence of two 
  time series}\label{sec:more}

The above test statistic can easily be adapted to test other
hypothesis. In this section, we consider one such
example, and test for
equality of serial dependence between two time series. Let us suppose
that $\{U_{t}\}$ and $\{V_{t}\}$ are two stationary time series, and
we wish to test whether they have the same sequential dependence
structure. Using the same motivation as that for the the goodness of
fit test described above we define the
test statistic 
\begin{eqnarray*}
\mathcal{P}_{T} = \frac{1}{T}\sum_{k=1}^{T}\int|\widehat{G}_{1,T}(x,y;\omega_{k}) - 
\widehat{G}_{2,T}(x,y;\omega_{k})|^{2}dF(x)dF(y),
\end{eqnarray*}
where $\widehat{G}_{1,T}$ and $\widehat{G}_{2,T}$ are the quantile
spectral density estimators based on $\{U_{t}\}$ and $\{V_{t}\}$
respectively and $F$ is any distribution function. 
In order to obtain the limiting distribution under the
null hypothesis we have 
$H_{0}: G_{1}(x,y;\omega) = G_{2}(x,y;\omega)$ and the alternative
 $H_{A}: G_{1}(x,y;\omega) \neq G_{2}(x,y;\omega)$ we expand
$\mathcal{P}_{T}$
\begin{eqnarray*}
\mathcal{P}_{T} := \mathcal{Q}_{1,1,T} +  \mathcal{Q}_{2,2,T} -
\mathcal{Q}_{1,2,T} - \mathcal{Q}_{2,1,T} + 2\mathcal{L}_{1,T} +
2\mathcal{L}_{2,T} + \mathcal{D},  
\end{eqnarray*} 
where 
\begin{eqnarray*}
\mathcal{Q}_{i,j,T} = 
\frac{1}{T}\sum_{k=1}^{T}\int\big[\widehat{G}_{i,T}(x,y;\omega_{k}) 
- 
\Ex(\widehat{G}_{i,T}(x,y;\omega_{k}))
\big]\big[\overline{\widehat{G}_{j,T}(x,y;\omega_{k}) - 
\Ex(\widehat{G}_{j,T}(x,y;\omega_{k})})
\big]dF_{}(x)dF_{}(y),
\end{eqnarray*}  
\begin{eqnarray*}
\mathcal{L}_{i,T} &=& \Re \frac{1}{T}\sum_{k=1}^{T}\int
\big[\widehat{G}_{i,T}(x,y;\omega_{k}) - 
\Ex(\widehat{G}_{i,T}(x,y;\omega_{k}))\big]\big[\Ex(\widehat{G}_{1,T}(x,y;\omega)) -
\Ex(\widehat{G}_{2,T}(x,y;\omega))\big]dF_{}(x)dF_{}(y)
\end{eqnarray*}
and 
\begin{eqnarray*}
\mathcal{D} = \int \int \int|\Ex(\widehat{G}_{1,T}(x,y;\omega)) - 
\Ex(\widehat{G}_{2,T}(x,y;\omega))|^{2}dF_{}(x)dF_{}(y)d\omega. 
\end{eqnarray*}
Therefore, using the above expansion
under the null hypothesis we have 
\begin{eqnarray*}
\mathcal{P}_{T} := \mathcal{Q}_{1,1,T} +  \mathcal{Q}_{2,2,T} -
\mathcal{Q}_{1,2,T} - \mathcal{Q}_{2,1,T},
\end{eqnarray*} 
where the moments are
%\begin{eqnarray*}
$\Ex\big(\mathcal{P}_{T}\big) 
= E_{T,3}+ O(\frac{1}{T}) = O(\frac{M}{T})$ and 
$\var\big(\mathcal{P}_{T} \big) 
= V_{T,3}    + O(\frac{1}{T})= O(\frac{M}{T^{2}})$,
%\end{eqnarray*}
with
\begin{eqnarray*}
E_{T,3} &=&  \frac{1}{T}\int \int W_{M}(\omega- \theta)^{2}\big(
G_{1}(x,x;\theta)G_{1}(y,y;\theta) +
G_{2}(x,x;\theta)G_{2}(y,y;\theta)\big)dF_{}(x)dF_{}(y)d\theta d\omega \\
V_{T,3} &=& \frac{4}{T^{2}}\sum_{i=1}^{2}\int \int
\Delta_{M}(\theta_{1} - \theta_{2})^2\prod_{j=1}^{2}
G_{i}(x_{1},y_{2};\theta_{i})G_{j}(y_{1},x_{2};\theta_{j})d\theta_{j}dF(x_{j})dF(y_{j}). 
\end{eqnarray*}

By using identical arguments as those used in the proof of Theorem
\ref{theorem:null}, under the null hypothesis we have 
\begin{eqnarray*}
V_{T,3}^{-1/2}\big(\mathcal{P}_{T} - E_{T,3}\big)\Dcon \mathcal{N}(0,1).
\end{eqnarray*}
Using the above result, we test for equality of sequential 
dependence, that is we reject the null hypothesis at the $\alpha$-level 
if $|V_{T,3}^{-1/2}(\mathcal{P}_{T} - E_{T,3})| > z_{1-\alpha}$. 

The limiting distribution of the alternative can be derived using 
the same methods as those used to derive the limiting distribution of 
$\mathcal{Q}_{T}$ under its alternative. It can be shown that 
\begin{eqnarray*}
\mathcal{P}_{T} - \mathcal{D} := \underbrace{2\mathcal{L}_{1,T} +
2\mathcal{L}_{2,T}}_{O_{p}(\frac{1}{\sqrt{T}})} 
+ O_{p}(\frac{M^{1/2}}{T}),
\end{eqnarray*} 
where $2\mathcal{L}_{1,T} +2\mathcal{L}_{2,T}$ can be approximated by
a quadratic form. Using this quadratic approximation, asymptotic 
normality of the above can be
shown. Thus under a fixed alternative the power grows to $100\%$ as
$T\rightarrow \infty$.

\begin{remark}
We can easily adapt our method to test that the distributions 
of $(X_{0},X_{r})$ and
$(X_{-r},X_{0})$ are identical for all $r$ (ie. $F_{r}(x,y) =
F_{-r}(x,y)$). This implies that the imaginary part of the quantile
spectral density $G(\cdot)$ is zero over all $x,y$ and $\omega$.
In this case, we use the test statistic
\begin{eqnarray*}
\mathcal{R}_{T} = \frac{1}{T}\sum_{r}\big|\Im \widehat{G}_{T}(x,y;\omega)\big|^{2}dF(x)dF(y),
\end{eqnarray*} 
where $F$ is some distribution, and by using identical methods to those
derived above we can obtain the limiting distribution of the test under the null. 
It is worth mentioning that \citeA{p:det-11} also
discuss the impact time reversibility has on the quantile
spectral density. 
\end{remark}

\section{Bootstrap approximation}\label{sec:practical}

The asymptotic normality result that we use to
obtain the p-value of the test statisic $\mathcal{Q}_{T}$ is only an
approximation. For small samples, the normality approximation may not be
particularly good, mainly because $\mathcal{Q}_{T}$ is a positive random
variable, whose distribution will be skewed. This may well
lead to more false positive than we can control for in our type I error. 

To correct for this, we propose estimating the finite sample distribution of
$\mathcal{Q}_{T}$ using a frequency domain bootstrap procedure. In a multivariate time
 series, the periodogram matrix at the 
fundamental frequencies asymptotically follow a Wishart distribution, 
moreover for our purposes they are close enough to be independent such that
we don't loose too much information by treating them as
independent (observe that the asymptotic variance of the test statistic
$\mathcal{Q}_{T}$ is only in terms of the pair-wise distributions and
does not contain any higher order dependencies). Thus motivated by the 
frequency domain bootstrap methods proposed in \citeA{p:hur-87}
and \citeA{p:fra-har-92} for univariate data and \citeA{p:ber-die-98} 
and \citeA{p:det-pap-09} for multivariate data, we propose the
following bootstrap scheme to obtain an estimate of the finite sample
distribution under the null hypothesis. 

Let $x_1 < \dots < x_q$ be a finite discretisation of the real line (noting that we
approximate $\mathcal{Q}_{T}$ with the discretisation 
\begin{eqnarray*}
\mathcal{Q}_{T} &=&
\frac{2\pi}{T} \sum_{k=1}^T \sum_{i_1, i_2=2}^q |\hat G(x_{i_1},
x_{i_2}; \omega_k) - \sum_{s}K_{M}(\omega_{k}-\omega_{s})G_{0}(x_{i,1},
x_{i_2}; \omega_k) |^2\times \\
&&(F_0(x_{i_1}) - F_0(x_{i_1-1}))( F_0(x_{i_2}) -
F_0(x_{i_2-1})).
\end{eqnarray*}
We observe that under the null hypothesis that 
$\bG_\bZ(\omega)$ will be the spectral density matrix of the
multivariate $q$-dimensional time series $\bZ_t  = (\tilde{Z}_t(x_1),
\cdots, \tilde{Z}_t(x_q))$ where $\tilde{Z}_t(x) = I(X_t \leq x) - F(x)$ and 
$\bG_\bZ(\omega)_{i_{1},i_{2}} = G_{0}(x_{i_{1}},x_{i_{2}};\omega)$. Thus
we use the transformation of $X_{t}$ into a high dimensional
multivariate time series to construct the the bootstrap distribution.

The steps of the frequency domain bootrap for the test statistic
$\mathcal{Q}_{T}$ are as follows:
\begin{itemize}
\item[] Step 1: Generate $T$ independent matrices 
$\bI^{*}_{\bZ}(\omega_k)  = \bG_{\bZ}(\omega_k)^{1/2 }W_{k}^{*}\bG_{\bZ}(\omega_k)
^{1/2 }$, where 
\begin{eqnarray*}
W_{k}^{*} \sim
\left\{
\begin{array}{cc}
W_{q}^{C}(1,I_{q}) &  1\leq k \leq T/2 \\
W_{q}^{R}(1,I_{q}) & k\in \{0,T/2\} \\
\overline{W_{T - k}^{*}} &  T/2 < k \leq T \\
\end{array}
\right.,
\end{eqnarray*}
and $W^{C}$ and $W^{R}$ denote the complex and real Wishart
distributions. 
\item[]Step 2: Construct the bootstrap quantile spectral density
  matrix estimators with 
$\hat \bG^{*}_\bZ(\omega_k) = \sum_s K_M(\omega_k - \omega_s)
\bI^{*}_{\bZ}(\omega_s) $ for $k=1, \dots, T$.

\item[]Step 3: 
Obtain the bootstrap test statistic 
\begin{eqnarray*}
Q_T^{*} = \frac{2\pi}{T} \sum_{k=1}^T \sum_{i_1, i_2=2}^q |\hat G^*(x_{i_1},
x_{i_2}; \omega_k) - G_0^M(x_{i_1},
x_{i_2}; \omega_k) |^{2}(F_0(x_{i_1}) - F_0(x_{i_1-1}))(F_0(x_{i_2}) -
F_0(x_{i_2-1})),
\end{eqnarray*} 
where $G_{0}^{M}(x,y;\omega) =
\frac{1}{2\pi}\sum_{k}\lambda(\frac{k}{M})C_{0,r}(x,y)\exp(ir\omega)$. 

\item[]Step 4: 
Approximate the distribution of $\mathcal{Q}_{T}$ under
  the null by using the empirical distribution of the bootstrap sample 
$\{\mathcal{Q}_{T}^{*}\}$. 

\item[]Step 5: Based on the bootstrap distribution estimate the
  p-value of $\mathcal{Q}_{T}$. 
\end{itemize}

\begin{figure}[h!]
\begin{center}
\includegraphics[width=14cm, height=5cm]{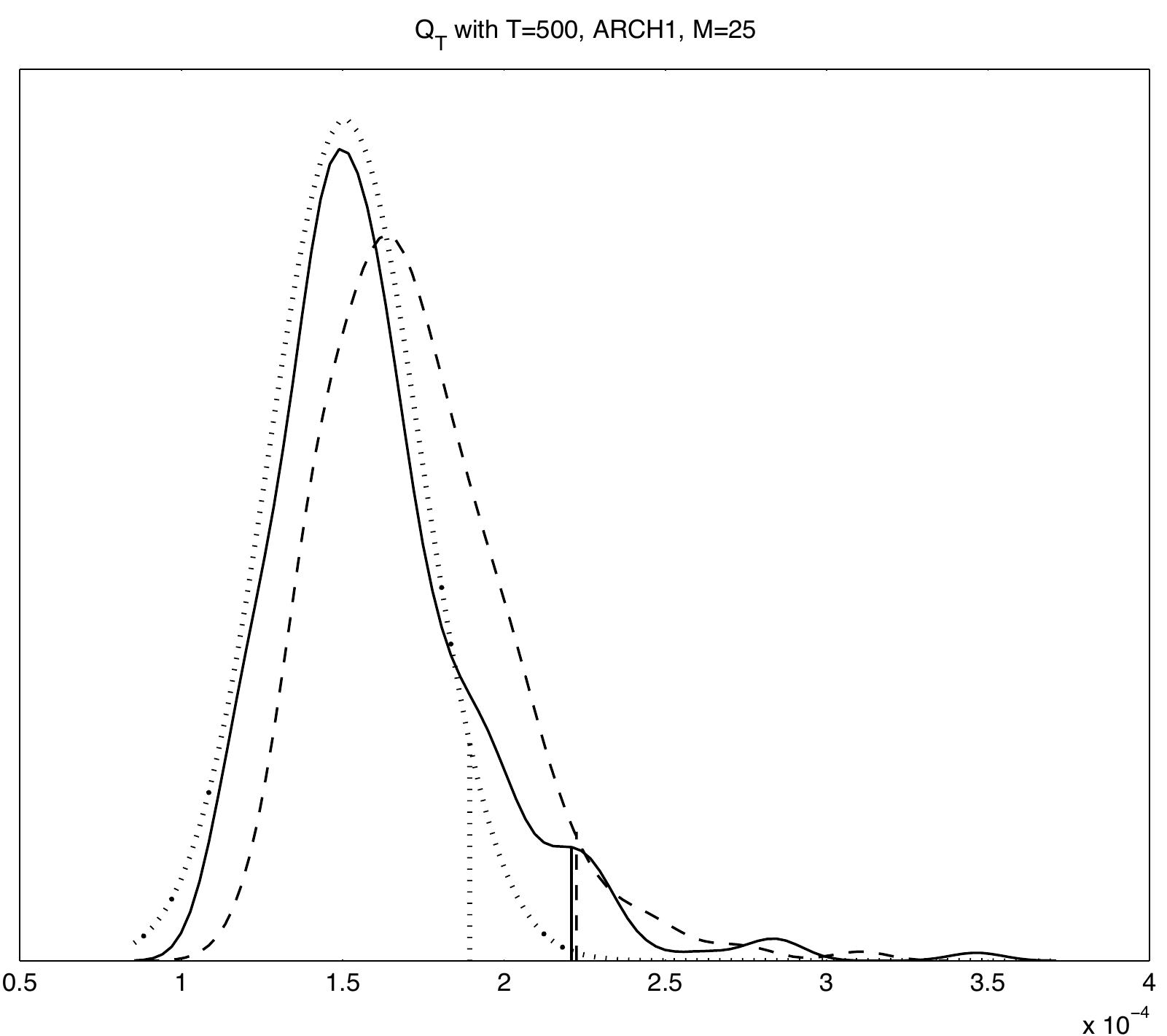}
\end{center}
\caption{\label{fig:6} The fine line is the standard normal (with the
  $5\%$ rejection line), the thick solid line is the finite sample
  density of the test statistic (with $5\%$ rejection region) and the
  thick dashed line is the bootstrap 
approximation (with $5\%$ rejection region).}
\end{figure}

We illustrate our procedure in Figure \ref{fig:6}, for this example we use
the quantile spectral density $G_{0}$, based on an ARCH$(1)$
($X_{t}=Z_{t}\sigma_{t}$ and $\sigma_{t}^{2} = a_{0}+a_{1}X_{t-1}^{2}$), where
$a_{0} = 1/1.9$, $a_{1} = 0.9$, $Z_{t}$ are iid standard normal random variables and $T=500$. A plot of the
normal approximation, the density of
$\mathcal{Q}_{T}$ (which is estimated and based on $500$ replications) 
and the bootstrap
estimator of the density (along with their rejection regions) 
is given in Figure \ref{fig:6}. 
We observe that the 
skew in the finite sample distribution means that the normal distribution 
is under estimating the location of the rejection region. However, 
the bootstrap approximation appears to capture relatively well the finite sample
distribution, and approximate well the rejection region.  
Since the bootstrap  scheme is based on sampling from iid random variables, 
we can write the bootstrap test statistic as a quadratic form. Thus by using
\citeA{p:lee-sub-11}, asymptotic normality of $\mathcal{Q}_{T}^{*}$ can be
shown with mean and variance given in (\ref{eq:ETstar}). Hence the limiting
distribution of the bootstrap statistic and limiting distribution of the test
statistic $\mathcal{Q}_{T}$, under the null, coincide. 

\section{Simulations and Real data examples}\label{sec:real}

\subsection{Simulations}

In this section we conduct a simulation study. In order to determine the 
effectiveness of the test we will use two different models that have the 
same first and second order structure (thus a test based on the covariance
structure would not be able to distinguish between them). 
In particular, we will consider 
the AR$(1)$ model $X_{t} = \mu + aX_{t-1} + \varepsilon_{t}$ and the 
squares of the ARCH$(1)$ model $Y_{t} = a_{0} + a Y_{t-1} + 
(Z_{t}^{2}-1)(a_{0} +aY_{t-1})$, where $\{\varepsilon_{t}\}$
and $\{Z_{t}\}$ are iid zero mean Gaussian random variables with 
$\var(Z_{t})=1$ and $\mu$ and $\var(\varepsilon_{t})$ chosen such that 
$X_{t}$ and $Y_{t}$ have the same mean and covariance structure. Note that 
in the simulation we only consider $a \leq 0.55$, so that the spectral
density of the squared ARCH exists. For each model we did 1000 replications
and the tests was done at both the $\alpha=0.1$ and $\alpha=0.05$ level. 

In our simulations we used the Bartlett window, compared the test for various
$M$ and used both the normal approximation and the proposed bootstrap
procedure. The results for $H_{0}:$ AR$(1)$ against the alternative $H_{A}:$
ARCH$(1)$ (various $a$, fixing $a_0 = 0.4$) are given in Table 
\ref{tab:1} and \ref{tab:2}. The results for 
$H_{0}:$ ARCH$(1)$ against $H_{A}:$ AR$(1)$ are given in Table
\ref{tab:3}  and \ref{tab:4}. 
We use the sample sizes $T=100$ and $500$. 

As expected under the null hypothesis the null hypothesis tends to over
reject, whereas the bootstrap gives a better approximation of the significance
level. There appears to be very little difference in the behaviour under the
null for various values of $a$ and between the AR and the ARCH. Under the 
alternative, the power seems to be quite high even for quite small samples. 
The only model where the power is not close to 100\% is when $a = 0.3$, sample
size $T=100$, the null is an 
AR$(1)$ and the alternative is an ARCH$(1)$. This can be explained by the fact
that for small values of $a$, both the AR and the ARCH models are 
relatively close to independent observations, thus making it relatively 
difficult to reject the null.  
%For $T=100$, the test statistics $\mathcal{Q}_T$ were evaluated using
%19 percentiles at $(5\%, 10\%, \dots, 95\%)$, whereas 49
%points at $(2\%, 4\%, \dots, 98\%)$ were used for $T=500$. For $H_0:
%ARCH(1)$, the true
%$p\%$ percentile was obtained by averaing 100 sample $p\%$ percentiles
%of simulated data with the sample size 8,000, and this method was also
%applied to estimate $G_0^M$. 

\subsection{Real Data}

%A plot of estimates of the quantile spectral density corresponding to the
%MSDR data is given in Figure \ref{fig:}

In this section we consider the the Microsoft daily return data (March, 1986 -
June, 2003) discussed in Section \ref{sec:motivate} and the Intel monthly
log return data (January 1973 - December 2003).
In the analysis below we will test whether the
GARCH and ARCH models are appropriate for the Microsoft and Intel data,
respectively. We use the Bartlett window. 

\begin{figure}[h!]
\begin{center}
\includegraphics[width=14cm, height=8cm]{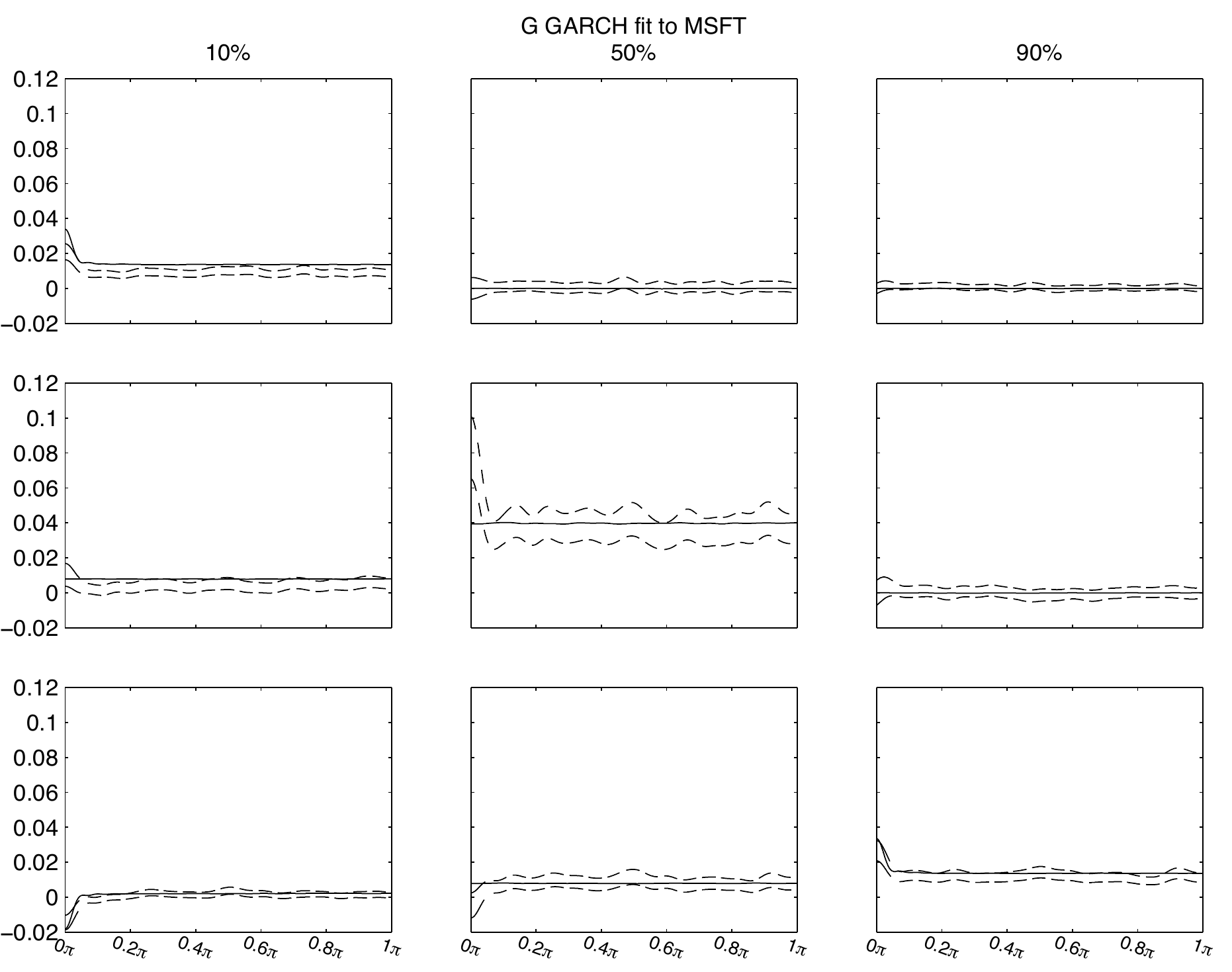}
\end{center}
\caption{\label{fig:7} The quantile spectral density of the fitted
  GARCH(1, 1) model using Microsoft data with the confidence intervals}
\end{figure}

A plot of the estimated $\widehat{G}_{T}$ together
with the piece-wise confidence intervals (obtained using the results
in Theorem \ref{theorem:QS}) and the corresponding quantile spectral 
density of the GARCH(1,1) is given in Figure \ref{fig:7} for the Microsoft
data. It is clear from the plot that
the GARCH$(1,1)$ model with coefficients evaluated using the maximum
likelihood estimator is not the appropriate model to fit
to this data. The plots suggest that the main deviation from the
GARCH$(1,1)$ arises at about $x,y=0$, 
This observation is confirmed by the results of
our test. Using various values of $M$ ranging from $30-70$, the
p-value corresponding to $\mathcal{Q}_{T}$ is almost zero with both the
normal approximation and also the Bootstrap method. Therefore,
from our analysis it seems that the GARCH$(1,1)$ is not a suitable
model for modelling the Microsoft daily returns from 1986-2003.

\begin{figure}[h!]
\begin{center}
\includegraphics[width=14cm, height=8cm]{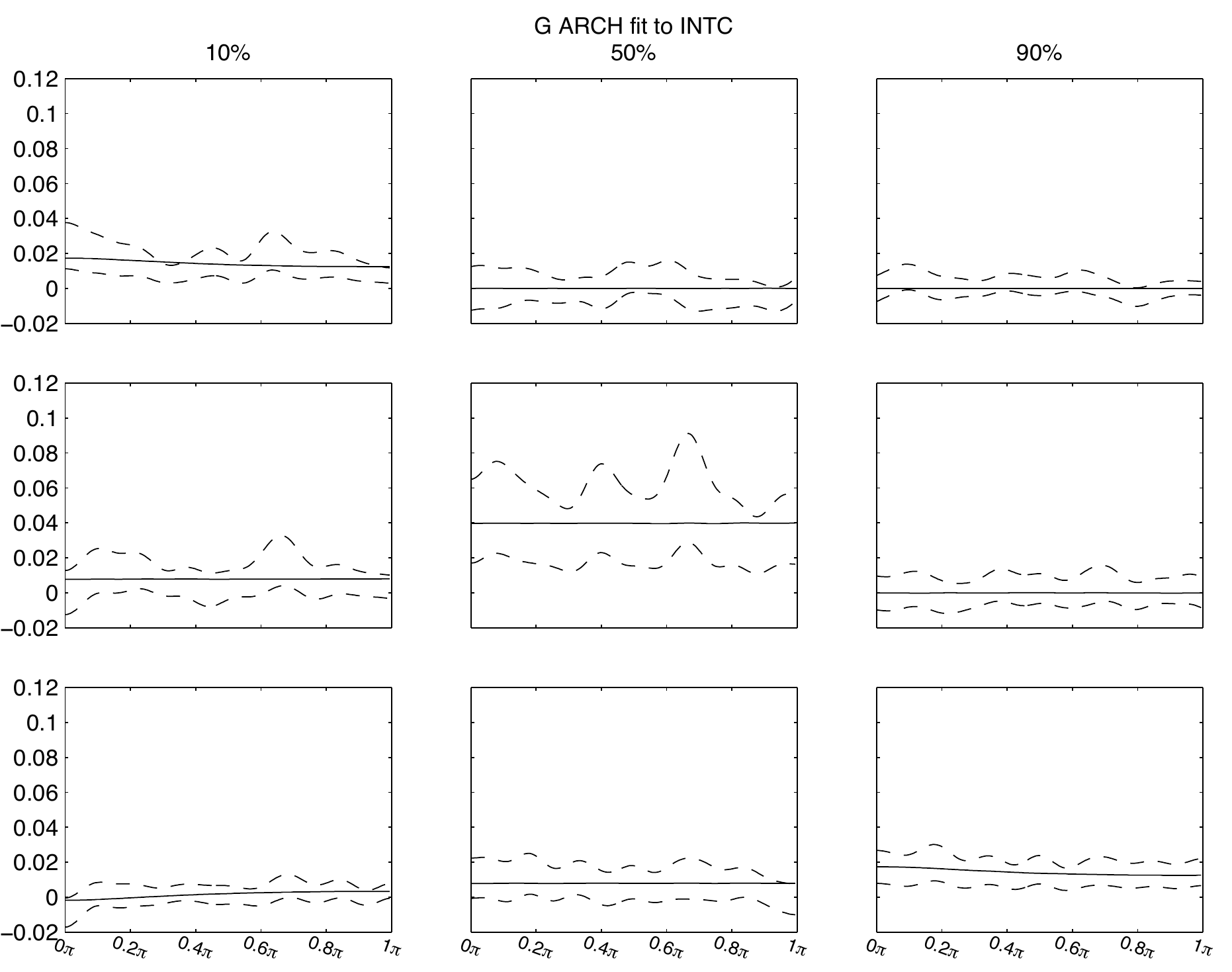}
\end{center}
\caption{\label{fig:8} The quantile spectral density of the fitted
  ARCH(1) from Intel data with the confidence intervals}
\end{figure}

We now consider the second data set, the Intel monthly log returns from 1973 -
2003. \citeA{b:tsa-05}
propose fitting an ARCH$(1)$ (with Gaussian innovations)
 model to this data, and maximum
likelihood yields the estimators $\mu = 0.0166$, $a_{0} = 0.0125$ and
$a_{1} = 0.363$, where $X_{t} =  \mu + \varepsilon_{t}$,
$\varepsilon_{t} = \sigma_{t}Z_{t}$ and $\sigma_{t}^{2} = a_{0} +
a_{1}\varepsilon_{t-1}^{2}$. A plot of the estimated $\widehat{G}_{T}$ 
with the piece-wise confidence intervals together  the quantile
spectral density of the ARCH($1$) model is
given in Figure \ref{fig:8}. We observe that  
quantile spectral density of the ARCH model 
lies in the confidence intervals for almost all
frequencies.  These observations are confirmed by 
the proposed goodness of fit test. A summary of
the results for various $M$, using both the normal approximation and
the bootstrap method is given in Table \ref{tab:5}. 
The p-values for the normal approximation tend to be smaller than the
p-values of the bootstrap method, this is probably due to the skew in
the finite sample distribution which results in smaller p-values. However, both the normal approximation
and the bootstrap give relatively large p-values for all values of
$M$.  Therefore there is not
enough evidence to reject the null. This backs the claims in \citeA{b:tsa-05}
that the ARCH$(1)$ may be an appropriate model for the the Intel
data. 

\begin{table}[h!]\small
\begin{center}
\begin{tabular}{c!{\VRule[0.8pt]}c!{\VRule[0.4pt]}c!{\VRule[0.4pt]}c!{\VRule[0.4pt]}c}
\specialrule{0.8pt}{0pt}{0pt}
   $M$    &          15    &         20    & 25    &         30 \\
\specialrule{0.8pt}{0pt}{0pt}
Normal 
p-value          & 0.0905 &  0.1279 &  0.1807 & 0.2643 \\ 
\hline
Bootstrap 
p-value &        0.3880 &     0.4320 &       0.4020 &
0.4780 \\
\specialrule{0.8pt}{0pt}{0pt}
\end{tabular}
\end{center}
\caption{The p-values for the Intel Data and various  values of $M$ \label{tab:5}}
\end{table}

% \begin{table}[h!]
% \begin{center}
% \begin{tabular}{|c|c|c|c|c|}
% \hline 
%    M    &          15    &         20    & 25    &         30 \\
% \hline
% \hline
% Normal 
% p-value          &0.0905 &  0.1279 &  0.1807 & 0.2643 \\ 
% \hline
% Bootstrap 
% p-value &        0.3880 &     0.4320 &       0.4020 &
% 0.4780 \\
% \hline
% \end{tabular}
% \end{center}
% \caption{The p-values for the Intel Data and various  values of $M$ \label{tab:5}}
% \end{table}

\subsection*{Acknowledgments}

This work has been partially supported by the National Science
Foundation, DMS-1106518.

\newpage

\appendix

\section{Proofs}

To obtain the sampling properties of $\hat{G}_{T}(\cdot)$ and $\mathcal{Q}_{T}$ (under both the
null and alternative), we first replace the
empirical distribution function $\hat{F}_{T}(x)$, with the
true distribution and show that the error is negligible. 
Define the zero mean, transformed variable 
$\tilde{Z}_{t}(x) = I(X_{t}\leq x) - F(x)$, where $F(\cdot)$ denotes
the marginal distribution of $\{X_{t}\}$. In addition define  
$\tilde{C}_{r}(x,y) =
\frac{1}{T}\sum_{t}\tilde{Z}_{t}(x)\tilde{Z}_{t+r}(y)$, 
\begin{eqnarray*}
\tilde{G}_{T}(x,y;\omega_{k}) &=&
\frac{1}{2\pi}\sum_{r}\lambda_{M}(r)\tilde{C}_{r}(x,y)\exp(ir\omega_{k}) =
\sum_{s}K_{M}(\omega_{k} - \omega_{s})\tilde{J}_{T}(x;\omega_{s})
\overline{\tilde{J}_{T}(y;\omega_{s})}, 
\end{eqnarray*}
\begin{eqnarray*}
\tilde{\mathcal{Q}}_{T} &=& 
\frac{1}{T}\sum_{k=1}^{T}\int|\tilde{G}_{T}(x,y;\omega_{k}) - 
\sum_{r}\lambda_{M}(r)C_{0,r}(x,y)\exp(ir\omega_{k})|^{2}dF_{0}(x)dF_{0}(y). 
\end{eqnarray*}
where $\tilde{J}_{T}(x;\omega) = \frac{1}{\sqrt{2\pi T}}
\sum_{t=1}^{T}Z_{t}(x)\exp(it\omega)$. 

In the proofs below we shall use the notation $\|X\|_{r} =
(\Ex(|X|^{r}))^{1/r}$. We first show that replacing $\hat{F}_{T}(x)$ with $F(x)$ does not affect the 
asymptotic sampling properties of $G_{T}(\cdot)$ and
$\mathcal{Q}_{T}$. 
\begin{lemma}\label{lemma:tildeQ}
Suppose Assumption \ref{assum:A} holds. Then we have
\begin{eqnarray}
\label{eq:Greplace}
(\Ex|\widehat{G}_{T}(x,y;\omega) - \tilde{G}_{T}(x,y;\omega)\big|^{2})^{1/2} = O(\frac{M}{T})
\end{eqnarray}
and 
\begin{eqnarray}
\label{eq:Qreplace}
(\Ex|\mathcal{Q}_{T} - \tilde{\mathcal{Q}}_{T}\big|^{2})^{1/2} = O(\frac{1}{T}). 
\end{eqnarray}
\end{lemma}
PROOF. We first observe that 
\begin{eqnarray*}
&&J_{T}(x;\omega_{k})\overline{J_{T}(y;\omega_{k})}  - 
\tilde{J}_{T}(x;\omega_{k})\overline{\tilde{J}_{T}(y;\omega_{k})}\\ 
&= & \left\{
\begin{array}{cc}
 0 & \omega_{k}\neq 0,\pi \\
 T  (\hat{F}_{T}(x) - F(x))(\hat{F}_{T}(y) - F(y)) & \textrm{
   otherwise }\\
\end{array}
\right.. 
\end{eqnarray*}
Substituting the above into $\hat{G}_{T}(\omega_{s}) -
\tilde{G}_{T}(\omega_{s})$ gives 
\begin{eqnarray}
\label{eq:Gdiff}
\hat{G}_{T}(\omega_{s}) - \tilde{G}_{T}(\omega_{s}) = 
TK_{M}(\omega_{s})(\hat{F}_{T}(x) - F(x))(\hat{F}_{T}(y) - F(y)). 
\end{eqnarray}
Using $K_M(\cdot) = O(\frac{M}{T})$ and $||\hat{F}_{T}(x)
- F(x) ||_2 = O(\frac{1}{T})$ in  (\ref{eq:Gdiff}), we obtain the
desired result for  (\ref{eq:Greplace}).  To prove (\ref{eq:Qreplace}) note that 
\begin{eqnarray*}
&&\mathcal{Q}_{T} - \tilde{\mathcal{Q}}_{T}  \\
&=& \hskip -2mm \int\frac{1}{T}\sum_{s=1}^{T}
\big(\hat{G}_{T}(x,y;\omega_{s}) - \tilde{G}_{T}(x,y;\omega_{s}) \big) 
\overline{\big(\hat{G}_{T}(x,y;\omega_{s}) + \tilde{G}_{T}(x,y;\omega_{s}) \big)}
dF_{0}(x)dF_{0}(y) \\ 
&+& \hskip -2mm  \Re\bigg(\int\frac{2}{T}\sum_{s=1}^{T}
\big(\hat{G}_{T}(x,y;\omega_{s}) - \tilde{G}_{T}(x,y;\omega_{s}) 
\big)G(x,y;\omega_{s}) 
dF_{0}(x)dF_{0}(y)\bigg).
\end{eqnarray*}
Thus substituting (\ref{eq:Gdiff}) into the above gives
\begin{eqnarray*}
&& \mathcal{Q}_{T} - \tilde{\mathcal{Q}}_{T} \\
&=& \int (\widehat{F}_{T}(x) - F(x))(\widehat{F}_{T}(y) - F(y)) \times\\ &&\bigg(\sum_{s=1}^{T}
K_{M}(\omega_{s})
\big(\hat{G}_{T}(x,y;\omega_{s}) + \tilde{G}_{T}(x,y;\omega_{s}) \big) \bigg)
dF_{0}(x)dF_{0}(y)  \\
 &+ & 2\int (\widehat{F}_{T}(x) - F(x))(\widehat{F}_{T}(y) - F(y))
\Re \bigg(\sum_{s=1}^{T}K_{M}(\omega_{s})G(x,y;\omega_{s}) \bigg)
dF_{0}(x)dF_{0}(y). 
\end{eqnarray*}
Therefore 
\begin{eqnarray*}
&&\big\|\mathcal{Q}_{T} - \tilde{\mathcal{Q}}_{T}\big\|_{2}\\
&\leq& \int  \big\|\widehat{F}_{T}(x) -
F(x)\big\|_{8}\big\|\widehat{F}_{T}(y) - F(y)\big\|_{8} \times \\ && 
\bigg(\sum_{s=1}^{T} \big( \big| K_{M}(\omega_{s})\big|\cdot
\big(\big\|\hat{G}_{T}(x,y;\omega_{s})\big\|_{8} + 
\|\tilde{G}_{T}(x,y;\omega_{s})\big\|_{8} \big) \big)
\bigg) dF_{0}(x)dF_{0}(y)  \\
&+&  2\int \big\|\widehat{F}_{T}(x) -
F(x)\big\|_{4}\big\|\widehat{F}_{T}(y) - F(y)\big\|_{4} \times \\ &&
\bigg(\sum_{s=1}^{T}\big|K_{M}(\omega_{s})|\cdot|G(x,y;\omega_{s})|\bigg)dF_{0}(x)dF_{0}(y). 
\end{eqnarray*}
For all $r\geq 2$, we have $\|\widehat{F}_{T}(x) - F(x)\|_{r} =
O(\frac{1}{\sqrt{T}})$, substituting this into the above gives 
$\big\|\mathcal{Q}_{T} -
\tilde{\mathcal{Q}}_{T}\big\|_{2}=O(\frac{1}{T})$, and the desired result. \hfill $\Box$

\vspace{3mm}
{\bf PROOF of Theorem \ref{theorem:QS}} To show asymptotic normality
of $\widehat{G}_{T}(\cdot)$, we first replace $\widehat{G}_{T}$ with $\tilde{G}_{T}$, by 
(\ref{eq:Greplace}) the replacement error is
$O_{p}(\frac{M}{T})$. Thus $\widehat{G}_{T}$ and $\tilde{G}_{T}$ have
the same asymptotic distribution and we can show how asymptotic
normality of $\widehat{G}_{T}$ by considering
$\tilde{G}_{T}(\cdot)$ instead. To show asymptotic normality of 
$\tilde{G}_{T}$ we use identical 
methods to those in Lee and Subba Rao (2011), where, since $\{I(X_{t} < x)\}$
are bounded random variables, we can use Ibragimov's covariance bounds for bounded random
variables. To obtain the limiting variance we note that
under Assumption \ref{assum:A}, since $s > 2$, we have that 
$\sum_{r}|r|\cdot|\cov(I(X_{0}\leq x),I(X_{r}\leq y))| < \infty$ and 
$\sum_{r_{1},r_{2},r_{2}}(1+|r_{j}|)|\cum(I(X_{0}\leq
x_{0}),I(X_{r_{1}}\leq x_{1}),I(X_{r_{2}}\leq x_{2}),
I(X_{r_{3}}\leq x_{3}))| < \infty$. Thus, the assumptions in
Brillinger (1981), Theorem 3.4.3 are satisfied, which allows us to
obtain the stated limiting variance. \hfill $\Box$

\vspace{3mm}
We use the following lemma to obtain a bound for the variance of
$\mathcal{Q}_{T}$. 

\begin{lemma}
Let the lag window be defined as in Definition \ref{def:lag} and
suppose $h_{1}(\cdot)$ and $h_{2}(\cdot)$ are bounded functions. Then we have 
\begin{eqnarray}
\label{eq:int1}
L_{1} = \int h_{1}(u_{1})h_{2}(u_{2})\Delta_{M}(u_{1}-u_{2})^{2}du_{1}du_{2} = 
O(M) 
\end{eqnarray}
and 
\begin{eqnarray}
\label{eq:int2}
L_{2} =\int h_{1}(u_{1})h_{2}(u_{2})\Delta_{M}(u_{1}+u_{2})\Delta_{M}(u_{1}-u_{2})
du_{1}du_{2} = O(1) 
\end{eqnarray}
where $\Delta_{M}(\cdot)$ is defined in (\ref{eq:delta}). 
\end{lemma}
PROOF. To simplify notation we prove the result for the truncated lag 
window $\lambda(u) = I_{[-1,1]}(u)$, but a similar result can also be proven 
for lag windows which satisfy Definition \ref{def:lag}. In the proof we
use the following two identities
\begin{eqnarray}
\sum_{t=0}^T e^{it\omega} =
e^{\frac{iT\omega}{2}} \frac{\sin( \frac{T+1}{2}  \omega)}{\sin
  (\omega/2)} \quad\textrm{ and }\quad \big(\int\big|\frac{\sin( \frac{M+1}{2}(u))}{
\sin((u)/2)}\big|^{p}du\big)^{1/p} = O(M^{1-p^{-1}}).\label{eq:fejer1}
\end{eqnarray}
We start by expanding $\Delta_{M}$ and using the above, to give 
\begin{eqnarray}
\Delta_{M}(\theta_{1} - \theta_{2}) &=& \int 
\sum_{j_{1},j_{2}=-M}^{M}\lambda_{M}(j_{1})\lambda_{M}(j_{2})
\exp(ij_{1}(\omega_{s_{1}}-\theta))\exp(ij_{2}(\omega_{s_{2}}-\theta))d\omega
\nonumber\\
&=& \sum_{j}\lambda_{M}(j)\lambda_{M}(-j)\exp(ij(\theta_{1}-\theta_{2})) 
 \nonumber\\
&=& \frac{\sin((M+1)(\theta_{1}-\theta_{2})/2)}{
\sin((\theta_{1}-\theta_{2})/2)} 2\Re e^{\frac{iM(\theta_{1}-\theta_{2})}{2}}. 
\label{eq:fejer}  
\end{eqnarray}
Substituting the above and (\ref{eq:fejer1}) into (\ref{eq:int1}) gives
\begin{eqnarray*}
|L_{1}| = \big|\int \int 
h_{1}(u_{1})h_{2}(u_{2})\Delta_{M}(u_{1}-u_{2})^{2}du_{1}du_{2}\big|
\leq \sup_{u,i}|h_{i}(u)|^{2}\int \int \big|\frac{\sin( \frac{M+1}{2}(u_{1}-u_{2}))}{
\sin((u_{1}-u_{2})/2)}\big|^{2}du_{1}du_{2} = O(M). 
\end{eqnarray*}
This proves (\ref{eq:int1}). To prove (\ref{eq:int2}) we observe that by a
change of variables ($v_{1} = u_{1} - u_{2}$ and $v_{2} = u_{1}+u_{2}$) we have 
\begin{eqnarray*}
|L_{2}| &\leq& C\int |\Delta_{M}(u_{1}+u_{2})|\cdot|\Delta_{M}(u_{1}-u_{2})|
du_{1} du_{2} \leq C\big(\int|\Delta_{M}(u)|du\big)^{2}.
\end{eqnarray*}
Now by substituting (\ref{eq:fejer}) and (\ref{eq:fejer1}) into the above 
gives $L_{2} = O(1)$. Thus we have obtained the desired result. \hfill $\Box$

\vspace{3mm}

{\bf PROOF Lemma \ref{lemma:mean-var}}
We first evalulate the expectation of $\mathcal{Q}_{T}$. By using Lemma \ref{lemma:tildeQ} we have  
\begin{eqnarray*}
&&\Ex(\mathcal{Q}_{T}) \\
&=& \frac{1}{T}\sum_{s=1}^{T}\int
\sum_{k_{1},k_{2}=1}^{T}K_{M}(\omega_{s}-\omega_{k_{1}})
K_{M}(\omega_{s}-\omega_{k_{2}})
\cov\big(\tilde{J}_{T}(x;\omega_{k_{1}})\overline{\tilde{J}_{T}(y;\omega_{k_{1}})}, 
\tilde{J}_{T}(x;\omega_{k_{2}})\overline{\tilde{J}_{T}(y;\omega_{k_{2}})}\big)  
\big) + O(\frac{1}{T}) \\
&=& I_{1} + I_{2} + I_{3} + O(\frac{1}{T}), 
\end{eqnarray*}
where 
\begin{eqnarray*}
I_{1} &=& \frac{1}{T}\int
\sum_{s,k_{1},k_{2}=1}^{T}\prod_{i=1}^{2}K_{M}(\omega_{s}- \omega_{k_{i}})
\cov\big(\tilde{J}_{T}(x;\omega_{k_{1}}), \tilde{J}_{T}(x;\omega_{k_{2}}))
\cov(\overline{\tilde{J}_{T}(y;\omega_{k_{2}})}, 
\overline{\tilde{J}_{T}(y;\omega_{k_{2}})})dF_{0}(x)dF_{0}(y)  \\
I_{2} &=& \frac{1}{T}\int
\sum_{s,k_{1},k_{2}=1}^{T}\prod_{i=1}^{2}K_{M}(\omega_{s}- \omega_{k_{i}})
\cov\big(\tilde{J}_{T}(x;\omega_{k_{1}}),\overline{\tilde{J}_{T}
(y;\omega_{k_{2}})})
\cov(\overline{\tilde{J}_{T}(y;\omega_{k_{1}})}, 
\tilde{J}_{T}(y;\omega_{k_{2}}))dF_{0}(x)dF_{0}(y) \\
I_{3} &=& \frac{1}{T}\int
\sum_{s,k_{1},k_{2}=1}^{T}\prod_{i=1}^{2}K_{M}(\omega_{s}- \omega_{k_{i}})
\cum\big(\tilde{J}_{T}(x;\omega_{k_{1}}),\overline{J}_{T}(y;\omega_{k_{1}}), 
\tilde{J}_{T}(x;\omega_{k_{2}}),\overline{\tilde{J}}_{T}(y;\omega_{k_{2}})\big)
dF_{0}(x)dF_{0}(y).  
\end{eqnarray*}
Under Assumption \ref{assum:A},  we have that 
$\sum_{r}|r|\cdot|\cov(I(X_{0}\leq x),I(X_{r}\leq y))| < \infty$ and 
$\sum_{r_{1},r_{2},r_{2}}(1+|r_{j}|)|\cum(I(X_{0}\leq
x_{0}),I(X_{r_{1}}\leq x_{1}),I(X_{r_{2}}\leq x_{2}),
I(X_{r_{3}}\leq x_{3}))| < \infty$. Therefore we can apply 
Brillinger (1981), Theorem 3.4.3 to obtain 
\begin{eqnarray*}
I_{1} &=&  \frac{1}{T}\sum_{s=1}^{T}\int
\sum_{k=1}^{T}K_{M}(\omega_{s}-\omega_{k})^{2}\int
G(x,x;\omega_{k})G(y,y;\omega_{k})dF_{0}(x)dF_{0}(y) + O(\frac{1}{T}) = O(\frac{M}{T}) \\
I_{2} &=& \frac{1}{T}\sum_{s=1}^{T}\int
\sum_{k=1}^{T}K_{M}(\omega_{s}-\omega_{k})K_{M}(\omega_{s}+\omega_{k})\int
G(x,y;\omega_{k})G(y,x;\omega_{k})dF_{0}(x)dF_{0}(y) + O(\frac{1}{T}) = O(\frac{1}{T}) \\
I_{3} &=& \frac{1}{T^{2}}\int\sum_{r}\lambda_{M}(r)^{2}
\sum_{t_{1},t_{2}=1}^{T}\cum(Z_{t_{1}}(x),Z_{t_{1}+r}(y),Z_{t_{2}}(x),
Z_{t_{2}+r}(y))dF_{0}(x)dF_{0}(y) = O(\frac{1}{T}).
\end{eqnarray*}
This gives us an asymptotic expression for the expectation. 
We now obtain an expression for the variance. Replacing $Z_{t}(\cdot)$
with $\tilde{Z}_{t}(\cdot)$ gives 
\begin{eqnarray*}
&&   \var(\mathcal{Q}_{T}) =   \\
&&   \frac{1}{T^2} \sum_{s_1,s_2=1}^T \int  \bigg(\sum_{k_1, k_2,
k_3, k_4} K_M(\omega_{s_1} - \omega_{k_1})   K_M(\omega_{s_1} - \omega_{k_2})
 K_M(\omega_{s_2} - \omega_{k_3})  K_M(\omega_{s_2} - \omega_{k_4})\\
&&   \times \cov\big( (J_{k_1,x_1} \overline{J}_{k_1,y_1} - 
\Ex(J_{k_1,x_1} 
\overline{J}_{k_1,y_1}))  (J_{k_2,x_1} 
\overline{J}_{k_2,y_1} - \Ex(J_{k_2,x_1} 
\overline{J}_{k_2,y_1})), \\ &&     (J_{k_3,x_2}
\overline{J}_{k_3,y_2} - \Ex(J_{k_3,x_2}
\overline{J}_{k_3,y_2}))(J_{k_4,x_2}
\overline{J}_{k_4,y_2} - \Ex(J_{k_4,x_2} 
\overline{J}_{k_4,y_2}))\big) \bigg)   dF_{0}(x_1) dF_{0}(y_1)dF_{0}(x_2) dF_{0}(y_2) 
+ O(\frac{1}{T}) \\
&&    =  II_{1} + II_{2} + II_{3} + O(\frac{1}{T})
\end{eqnarray*}
where $J_{k,x} = \tilde{J}_{T}(x;\omega_{k})$,
\begin{eqnarray*}
II_{1}  &=& \frac{1}{T^2} \sum_{s_1,s_2=1}^T \int\sum_{k_1,k_2,k_3,k_4} 
 \cum(J_{k_{1},x_{1}}\bar{J}_{k_{1},y_{1}}, \bar{J}_{k_{3},x_{2}}J_{k_{3},y_{2}})\cum(J_{k_{2},x_{1}}\bar{J}_{k_{2},y_{1}},\overline{J}_{k_{4},x_{2}}J_{k_{4},y_{2}}) \\
&&\prod_{i=1}^{2}K_{M}(\omega_{s_{1}}  - \omega_{k_{i}}) \prod_{i=3}^{4}K_{M}(\omega_{s_{2}}  - \omega_{k_{i}}) dF_{0}(x_1) dF_{0}(y_1)dF_{0}(x_2) dF_{0}(y_2) \\
II_{2} &=& 
\frac{1}{T^2} \sum_{s_1,s_2=1}^T \int\sum_{k_1,k_2,k_3,k_4} 
 \cum(J_{k_{1},x_{1}}\bar{J}_{k_{1},y_{1}},\overline{J}_{k_{4},x_{2}}J_{k_{4},y_{2}})\cum(J_{k_{2},x_{1}}\bar{J}_{k_{2},y_{1}},
 \bar{J}_{k_{3},x_{2}}J_{k_{3},y_{2}})\\ 
&&\prod_{i=1}^{2}K_{M}(\omega_{s_{1}}  - \omega_{k_{i}}) \prod_{i=3}^{4}K_{M}(\omega_{s_{2}}  - \omega_{k_{i}}) dF_{0}(x_1) dF_{0}(y_1)dF_{0}(x_2) dF_{0}(y_2) \\ 
II_{3}  &=&  \frac{1}{T^2} \sum_{s_1,s_2=1}^T \int\sum_{k_1,k_2,k_3,k_4} 
\cum(J_{k_{1},x_{1}}\bar{J}_{k_{1},y_{1}},
J_{k_{2},x_{1}}\bar{J}_{k_{2},y_{1}},\bar{J}_{k_{3},x_{2}}J_{k_{3},y_{2}},
\overline{J}_{k_{4},x_{2}}J_{k_{4},y_{2}}) \\
&& \prod_{i=1}^{2}K_{M}(\omega_{s_{1}}  - \omega_{k_{i}}) \prod_{i=3}^{4}K_{M}(\omega_{s_{2}}  - \omega_{k_{i}}) dF_{0}(x_1) dF_{0}(y_1)dF_{0}(x_2) dF_{0}(y_2).
\end{eqnarray*}
To obtain an expression for the variance we start by expanding $II_{1}$
\begin{eqnarray*}
&II_{1}& = \frac{1}{T^2} \sum_{s_1, s_2}\int \sum_{k_1,
  k_2, k_3, k_4}
\prod_{i=1}^{2}K_{M}(\omega_{s_{1}}  - \omega_{k_{i}})
\prod_{i=3}^{4}K_{M}(\omega_{s_{2}}  - \omega_{k_{i}})\\
& \times & \bigg(
\cov(J_{k_1,x_1} , J_{k_3,x_2}) \cov(\overline J_{k_1,y_1},\overline
J_{k_3,y_2}  ) \cov(J_{k_2,x_1} , J_{k_4,x_2})\cov(\overline J_{k_2,y_1},\overline J_{k_4,y_2}  )
\\
& +&   \cov(J_{k_1,x_1} , J_{k_3,x_2})
\cov(\overline J_{k_1,y_1},\overline J_{k_3,y_2}  ) \cov( J_{k_2,x_1} ,
\overline J_{k_4,y_2}) \cov(\overline J_{k_2,y_1} , J_{k_4,x_2} )
\\
& +& \cov(J_{k_1,x_1} , J_{k_3,x_2})
\cov(\overline J_{k_1,y_1},\overline J_{k_3,y_2}  ) \cum(J_{k_2,x_1} ,
\overline J_{k_2,y_1} ,J_{k_4,x_2}, \overline J_{k_4,y_2} )\\
& +& 
\cov( J_{k_1,x_1} ,  \overline J_{k_3,y_2}) \cov(\overline J_{k_1,y_1}
, J_{k_3,x_2} )  \cov(J_{k_2,x_1} , J_{k_4,x_2})\cov(\overline J_{k_2,y_1},\overline J_{k_4,y_2}  )
\\
& +&  \cov( J_{k_1,x_1} ,  \overline J_{k_3,y_2}) \cov(\overline J_{k_1,y_1} , J_{k_3,x_2} )  \cov( J_{k_2,x_1} ,  \overline
J_{k_4,y_2}) \cov(\overline J_{k_2,y_1} , J_{k_4,x_2} )
\\
& +& 
\cov( J_{k_1,x_1} ,  \overline
J_{k_3,y_2}) \cov(\overline J_{k_1,y_1} , J_{k_3,x_2} )  \cum(
J_{k_2,x_1} , \overline J_{k_2,y_1} ,
J_{k_4,x_2}, \overline J_{k_4,y_2} )\\ 
& +& 
\cum(
J_{k_1,x_1} , \overline J_{k_1,y_1} ,
J_{k_3,x_2}, \overline J_{k_3,y_2} )  \cov(J_{k_2,x_1} , J_{k_4,x_2})\cov(\overline
J_{k_2,y_1},\overline J_{k_4,y_2}  )
\\
& +&  \cum(
J_{k_1,x_1} , \overline J_{k_1,y_1} ,
J_{k_3,x_2}, \overline J_{k_3,y_2} )   \cov( J_{k_2,x_1} ,  \overline
J_{k_4,y_2}) \cov(\overline J_{k_2,y_1} , J_{k_4,x_2} )
\\
& +& 
\cum(
J_{k_1,x_1} , \overline J_{k_1,y_1} ,
J_{k_3,x_2}, \overline J_{k_3,y_2} )   \cum(
J_{k_2,x_1} , \overline J_{k_2,y_1} ,
J_{k_4,x_2}, \overline J_{k_4,y_2} )
\bigg) \prod_{j=1}^{2}dF_{0}(x_j)dF_{0}(y_j) \\ & := & \sum_{j=1}^{9}II_{1,j}.
\end{eqnarray*}
We use \cite{b:bri-81}, Theorem 3.4.3 to obtain the following expression for $II_{1,1}$
\begin{eqnarray*}
 II_{1,1} \hskip -3mm &=&\frac{1}{T^2} \sum_{s_1, s_2}\int \bigg( \sum_{k_1, k_{2}=1}^{T}
\big(\prod_{i=1}^2 K_M(\omega_{s_i} - \omega_{k_1})   K_M(\omega_{s_i} - \omega_{k_2})\big)\\
&&\hskip -9mm \cov(J_{k_1,x_1} , J_{k_1,x_2}) \cov(\overline J_{k_1,y_1},\overline
J_{k_1,y_2}) 
\cov(J_{k_2,x_1} , J_{k_2,x_2})\cov(\overline J_{k_2,y_1},\overline
J_{k_2,y_2}) \bigg)\\
&&\hskip -9mm\prod_{j=1}^{2}dF_{0}(x_j)dF_{0}(y_j)  + O(\frac{1}{T^2})
\\ 
=&& \hskip -9mm\frac{1}{T^2} \int\int \bigg( \int W_M(\omega_{s_1} - \theta_1)
W_M(\omega_{s_1} - \theta_2) d\omega_{s_1}  \bigg)\times\\
&&\hskip -9mm \bigg( \int W_M(\omega_{s_2} - \theta_1)
W_M(\omega_{s_2} - \theta_2) d\omega_{s_2}  \bigg)  \prod_{i=1}^{2}
G(x_{1},x_{2};\theta_{i})G(y_{1},y_{2}; -\theta_{i})d\theta_{i} \\ 
&&\hskip -9mm \prod_{j=1}^{2}dF_{0}(x_j)dF_{0}(y_j)  + O(\frac{1}{T^2}) \\
=& &\hskip -9mm\frac{1}{T^{2}}\int \int
\Delta_{M}(\theta_{1} - \theta_{2})^{2}\prod_{i=1}^{2}
G(x_{1},x_{2};\theta_{i})G(y_{1},y_{2}; -\theta_{i})d\theta_{i}
\prod_{j=1}^{2}dF_{0}(x_{j})dF_{0}(y_{j})+ O(\frac{1}{T^{2}}).
\end{eqnarray*} 
Therefore by using (\ref{eq:int1}) we have $II_{1,1}  = 
O(\frac{M}{T^{2}})$. We now consider $II_{1,2}$, by
using a similar method we have 
\begin{eqnarray*}
II_{1,2} &=&
\frac{1}{T^2}\int \bigg( W_M(\omega_{s_1} - \theta_1) W_M(\omega_{s_1} -
\theta_2) W_M(\omega_{s_2} - \theta_1) W_M(\omega_{s_2} + \theta_2) \times \\ 
&&G(x_1, x_2 ,\theta_1) G(y_1, y_2 , -\theta_1) G(x_1, y_2 ,
\theta_2) G(y_1, x_2 , -\theta_2) \bigg) d\theta_1 d\theta_2 d\omega_{s_1}
d\omega_{s_2}  
\prod_{j=1}^{2}dF_{0}(x_{j})dF_{0}(y_{j})  + O(\frac{1}{T^2})\\
 &=& \frac{1}{T^2}\int \Delta_{M}(\theta_{1}-\theta_{2})\Delta_{M}(\theta_{1}+\theta_{2})
G(x_1, x_2 ,\theta_1) G(y_1, y_2 , -\theta_1) \\ && G(x_1, y_2 ,
\theta_2) G(y_1, x_2 , -\theta_2) d\theta_1 d\theta_2 \prod_{j=1}^{2}dF_{0}(x_{j})dF_{0}(y_{j}) + O(\frac{1}{T^{2}}). 
\end{eqnarray*}
By using (\ref{eq:int2}) the above integral is $O(1)$, and 
altogether $II_{1,2} = O(\frac{1}{T^{2}})$. Using a similar argument,
one can show that $II_{1,3}$, $II_{1,4}$ are smaller than
$O(\frac{M}{T^2})$, so negligible. 
For $II_{1,5}$, we use that
\begin{eqnarray*}
 \cov(J_{k_1,x}, \overline{J}_{k_2,y})  =
 \begin{cases}
G(x, y, \omega_{k_1})  & k_1 + k_2 = T\\
O(\frac{1}{T}) & otherwise   
 \end{cases}
\end{eqnarray*}
, which follows from \cite{b:bri-81}, Theorem 3.4.3. This leads to 
\begin{eqnarray*}
 II_{1,5} &=&\frac{1}{T^2} \sum_{s_1, s_2}\int \bigg( \sum_{k_1, k_{2}=1}^{T}
\big(\prod_{i=1}^2 K_M(\omega_{s_1} - \omega_{k_i})   K_M(\omega_{s_2} + \omega_{k_i})\big)\\
&&\cov(J_{k_1,x_1} , J_{k_1,y_2}) \cov(\overline J_{k_1,y_1},\overline
J_{k_1,x_2}) 
\cov(J_{k_2,x_1} , J_{k_2,y_2})\cov(\overline J_{k_2,y_1},\overline
J_{k_2,x_2}) \bigg)\\
&&\prod_{j=1}^{2}dF_{0}(x_j)dF_{0}(y_j)  + O(\frac{1}{T^2})
\\
&= & \frac{1}{T^2} \int\int \bigg( \int W_M(\omega_{s_1} - \theta_1)
W_M(\omega_{s_1} - \theta_2) d\omega_{s_1}  \bigg)\times\\
&&\bigg( \int W_M(\omega_{s_2} + \theta_1)
W_M(\omega_{s_2} + \theta_2) d\omega_{s_2}  \bigg)  \prod_{i=1}^{2}
G(x_{1},y_{2};\theta_{i})G(y_{1},x_{2}; -\theta_{i})d\theta_{i} \\ 
&&\prod_{j=1}^{2}dF_{0}(x_j)dF_{0}(y_j)  +
O(\frac{1}{T^2}) \\
&=& II_{1,1}
\end{eqnarray*}
because of $\Delta(\theta) = \Delta(-\theta)$ and interchangeability
of integrals about $(x_1, x_2, y_1, y_2)$. 
With a similar method, one can show that $II_{1,6} \dots,II_{1,9}$ are all dominated by
$II_{1,1}$ and $II_{1,5}$ 
Altogether this gives 
\begin{eqnarray*}
II_{1} = \frac{2}{T^{2}}\int
\Delta_{M}(\theta_{1} - \theta_{2})^{2}\prod_{i=1}^{2}
G(x_{1},x_{2};\theta_{i})G(y_{1},y_{2};\theta_{i})d\theta_{i}
\prod_{j=1}^{2}dF_{0}(x_{j})dF_{0}(y_{j}) + O(\frac{1}{T^{2}}).
\end{eqnarray*}
Using the identical argument with the above,  we can show that
\begin{eqnarray*}
II_{2} =\frac{2}{T^{2}}\int
\Delta_{M}(\theta_{1} - \theta_{2})^{2}\prod_{i=1}^{2}
G(x_{1},x_{2};\theta_{i})G(y_{1},y_{2};\theta_{i})d\theta_{i}
\prod_{j=1}^{2}dF_{0}(x_{j})dF_{0}(y_{j})+ O(\frac{1}{T^{2}}).
\end{eqnarray*}
To bound $II_{3}$ we recall that 
\begin{eqnarray*}
II_{3} &=& \frac{1}{T^2} \sum_{s_1,s_2=1}^T \int  \sum_{k_1, k_2,
k_3, k_4} K_M(\omega_{s_1} - \omega_{k_1})   K_M(\omega_{s_1} - \omega_{k_2})
 K_M(\omega_{s_2} - \omega_{k_3})  K_M(\omega_{s_2} - \omega_{k_4})\\
 && \cum\big(J_{k_1,x_1} \overline{J}_{k_1,y_1},  
J_{k_2,x_1}\overline{J}_{k_2,y_1},J_{k_3,x_2}
\overline{J}_{k_3,y_2},J_{k_4,x_2}\overline{J}_{k_4,y_2}\big) 
dF_{0}(x_1) dF_{0}(y_1)dF_{0}(x_2) dF_{0}(y_2) + O(\frac{1}{T}). 
\end{eqnarray*}
By using the method of indecomposable partitions (see \cite{b:bri-81}, 
Theorem 2.3.2) to partition the above cumulant of products into the
product of cumulants. 
This together with \cite{b:bri-81}, Theorem 3.4.3 gives us
$II_{3} = O(\frac{M}{T^{3}})$ (see Lee's thesis for further details).

Combining the expressions for $II_{1}$, $II_{2}$ and $II_{3}$ gives us 
the expression for the variance and completes the proof. \hfill $\Box$

\vspace{3mm}

\subsection{Proof of Theorem \ref{theorem:null}}

Now we show that $\mathcal{Q}_{T}$ can be approximated by the sum of
martingale differences, this will allow us to use the martingale
central limit theorem to prove Theorem \ref{theorem:null}.

We first define the martingale difference decomposition of
$\tilde{Z}_{t}(x) = \sum_{j=0}^{t}M_{j}^{(x)}(t-j)$, where
$M_{j}^{(x)}(t-j) = \Ex(Z_{t}(x)|\mathcal{F}_{t-j}) -
\Ex(Z_{t}(x)|\mathcal{F}_{t-j-1})$, where for $t > 0$ we have 
$\mathcal{F}_{t} =
\sigma(X_{t},X_{t-1},\ldots,X_{1})$  and for $t \leq 0$ we let
$\mathcal{F}_{t} = \sigma(1)$, and $M_{j}(s) = 0$ for $j \geq s$. 
Using the above notation we define the random variable   
\begin{eqnarray}
\label{eq:mathcalS}
\mathcal{S}_{T} &=& \frac{1}{T^{2}}\int \sum_{j_{1},\ldots,j_{4}=0}^{\infty}
\sum_{t_{1},r,t_{2}\in \mathcal{A}}\lambda_{M}(r)^{2}
M^{(x)}_{j_{1}}(t_{1}-j_{1})M_{j_{2}}^{(y)}(t_{1}+r-j_{2})\nonumber\\
&&\times
M^{(x)}_{j_{3}}(t_{2}-j_{3})M^{(y)}_{j_{4}}(t_{2}+r-j_{4})dF_{0}(x)dF_{0}(y),
\end{eqnarray}
where $\mathcal{A}$ = $\{(t_{1}-j_{1},t_{1}+r-j_{2},t_{2}-j_{3},
t_{2}+r-j_{4})$ are all different $\}$.  
\begin{theorem}\label{thm:bounds}
Suppose Assumption \ref{assum:A} holds, $\mathcal{S}_{T}$ is
defined as in (\ref{eq:mathcalS}) and the null hypothesis is true. Then we have 
\begin{eqnarray*}
\mathcal{Q}_{T} - \Ex(\mathcal{Q}_{T}) &=& \mathcal{S}_{T} +
O_{p}(\frac{1}{T}+\frac{M^{1/2}}{T^{3/2}})
\end{eqnarray*} 
and for all $r \geq 2$
$\big\|\mathcal{S}_{T}\big\|_{r} = O(\frac{M^{1/2}}{T})$. 
\end{theorem}
PROOF. We prove the result using a combination of iterative martingales and
Burkholder's inequality for martingale differences. 
We note that for $r \geq 2$ we have 
\begin{eqnarray}
\|M_{j}^{(x)}(t-j)\|_{r} &=& \|\Ex(Z_{t}(x)|\mathcal{F}_{t-j}) - 
\Ex(Z_{t}(x)|\mathcal{F}_{t-j-1})\|_{2}\leq 
2\|\Ex(Z_{t}(x)|\mathcal{F}_{t-j})\|_{r}
\leq C\alpha(j), \label{eq:Mjb} 
\end{eqnarray}
where $\mathcal{F}_{t}=\sigma(X_{t},X_{t-1},\ldots,X_{1})$, which follows from 
Ibragimov's inequality. 
Substituting the representation $Z_{t}(x)  = \sum_{j=1}^{\infty}M_{j}^{(x)}(t-j)$
into $\mathcal{Q}_{T}$ gives 
\begin{eqnarray*}
&&\mathcal{Q}_{T} - 
\Ex(\mathcal{Q}_{T}) \\
&=& \frac{1}{T^{2}}\sum_{j_{1},\ldots,j_{4}=0}^{\infty}
\sum_{r=-M}^{M}\lambda_{M}(r)^{2}\sum_{t_{1},t_{2}}\bigg(
\overline{M^{(x)}_{j_{1}}(t_{1}-j_{1})M_{j_{2}}^{(y)}(t_{1}+r-j_{2})} \times
\overline{M^{(x)}_{j_{3}}(t_{2}-j_{3})M_{j_{4}}^{(y)}(t_{2}+r-j_{4})}\\
&-&  \Ex\big(\overline{M^{(x)}_{j_{1}}(t_{1}-j_{1})M_{j_{2}}^{(y)}(t_{1}+r-j_{2})} 
\times
\overline{M^{(x)}_{j_{3}}(t_{2}-j_{3})M_{j_{4}}^{(y)}(t_{2}+r-j_{4})})\big)
\bigg)dF_{0}(x)dF_{0}(y),
\end{eqnarray*}
where $\overline{X}$ denotes the centralised random variable
$X-\Ex(X)$ (note that $M_{j}(s) = 0$ for $s \leq 0$). We now partition the above sum into several cases, where we treat 
$j_{1},\ldots,j_{4}$ as free and condition on $t_{1},t_{2}$ and $r$:
\begin{itemize}
\item[(i)] $\mathcal{A}$ = 
$\{(t_{1},t_{2},r) \textrm{ such that } (t_{1}-j_{1},t_{1}+r-j_{2},t_{2}-j_{3},t_{2}+r-j_{4})$ are all different $\}$.
\item[(ii)] $\mathcal{B}$ = $\{(t_{1},t_{2},r) \textrm{ such that }(t_{1}-j_{1}=t_{1}+r-j_{2})$ and 
$(t_{2}-j_{3}=t_{2}+r-j_{4})\}$.
\item[(iii)] $\mathcal{C} = \{(t_{1},t_{2},r) \textrm{ such that }(t_{1}-j_{1}) = (t_{2}-j_{3})$ or
  $(t_{2}+r-j_{4})$ and $(t_{1}+r-j_{2})\neq (t_{1}-j_{1})\}$. 
\item[(iv)] $\mathcal{D} = \{(t_{1},t_{2},r) \textrm{ such that }(t_{1}+r-j_{2}) = (t_{2}-j_{3})$ or
  $(t_{2}+r-j_{4})$ and $(t_{1}-j_{1})\neq (t_{1}+r-j_{2})\}$.
\item[(v)]$ \mathcal{E} = \{(t_{1},t_{2},r) \textrm{ such that }(t_{2}-j_{3}) = (t_{1}-j_{1})$ or
  $(t_{1}+r-j_{2})$ and $(t_{2}+r-j_{4}\neq t_{2}-j_{3})\}$. 
\item[(iv)] $\mathcal{F} = \{(t_{1},t_{2},r) \textrm{ such that }(t_{2}+r-j_{4}) = (t_{1}-j_{2})$ or
  $(t_{1}+r-j_{2})$ and $(t_{2}-j_{3})\neq (t_{2}+r-j_{4})\}$.  
\end{itemize}
%\begin{eqnarray*}
%\mathcal{C} = 
%\left\{
%t_{1}-j_{1}= 
%\left\{
%\begin{array}{c}
%t_{2}-j_{3} \\
%\textrm{or} \\
%t_{2}+r-j_{4} \\
%\end{array}
%\right.
%\right.
%\end{eqnarray*}
%\item[(iv)] 
%\begin{eqnarray*}
%\mathcal{D} = 
%\left\{t_{1}+r-j_{2} = 
%\left\{
%\begin{array}{c}
%t_{2}-j_{3} \\
%\textrm{or} \\
%t_{2}+r-j_{4} \\
%\end{array}
%\right.
%\right.
%\end{eqnarray*}
%\item[(v)] 
%\begin{eqnarray*}
%\mathcal{E} &=& 
%\left\{
%\begin{array}{ll}
%t_{1}-j_{1} = t_{1}+r-j_{2} = t_{2}-j_{3} & t_{2}+r-j_{4}\textrm{
%  different}\\
%t_{1}-j_{1} = t_{1}+r-j_{2} = t_{2}+r-j_{4} & t_{2}-j_{3}\textrm{
%  different} \\
%t_{1}-j_{1} = t_{2}-j_{3}  = t_{2}+r - j_{4} & t_{1}+r-j_{2}\textrm{
%  different} \\
%t_{1}+r-j_{2} = t_{2}-j_{3} = t_{2}+r-j_{4} & t_{1}-j_{2}\textrm{
%  different}
%\end{array}
%\right.
%\end{eqnarray*}
%\item[(vi)] $\mathcal{F} = \{t_{1}-j_{1} = t_{1}+r-j_{2} = t_{2}-j_{3} 
%= t_{2}+r-j_{4}\}$.  
%\end{itemize}
Thus 
\begin{eqnarray*}
\mathcal{Q}_{T} - \Ex(\mathcal{Q}_{T}) = 
\int \bigg(I_{\mathcal{A}} +
I_{\mathcal{B}} + I_{\mathcal{C}} + I_{\mathcal{D}} + I_{\mathcal{E}} +I_{\mathcal{F}}\bigg)
dF_{0}(x)dF_{0}(y),     
\end{eqnarray*}
where  
\begin{eqnarray*}
I_{\mathcal{A}} &=& \frac{1}{T^{2}}\sum_{j_{1},\ldots,j_{4}=0}^{\infty}
\sum_{r,t_{1},t_{2}\in \mathcal{A}}\lambda_{M}(r)^{2}
M^{(x)}_{j_{1}}(t_{1}-j_{1})M_{j_{2}}^{(y)}(t_{1}+r-j_{2})
M^{(x)}_{j_{3}}(t_{2}-j_{3})M^{(y)}_{j_{4}}(t_{2}+r-j_{4}), 
\end{eqnarray*}
\begin{eqnarray*}
I_{\mathcal{B}} &=&  \frac{1}{T^{2}}\sum_{j_{1},\ldots,j_{4}=0}^{\infty}
\sum_{r,t_{1},t_{2}\in \mathcal{B}}\lambda_{M}(j_{1}-j_{2})^{2} \bigg(
\overline{M^{(x)}_{j_{1}}(t_{1}-j_{1})M_{j_{2}}^{(y)}(t_{1}-j_{1})}
\times 
\overline{M^{(x)}_{j_{3}}(t_{2}-j_{3})M_{j_{4}}^{(y)}(t_{2}-j_{3})} -  \\
&&\Ex\big(
\overline{M^{(x)}_{j_{1}}(t_{1}-j_{1})M_{j_{2}}^{(y)}(t_{1}-j_{1})} \times
\overline{M^{(x)}_{j_{3}}(t_{2}-j_{3})M_{j_{4}}^{(y)}(t_{2}-j_{3})}
\big) \bigg) 
\end{eqnarray*}
for $I_{\mathcal{C}},\ldots,I_{\mathcal{F}}$ are defined similarly. 

We first bound $I_{\mathcal{A}}$. We partition $\mathcal{A}$ into $24$
cases by the order of $
(t_{1}-j_{1},t_{1}+r-j_{2},t_{2}-j_{3},t_{2}+r-j_{4})$. The first is $\mathcal{A}_{1} = \{(t_{1},t_{2},r) \textrm{ such that }t_{1} - j_{1} > t_{1}+r-j_{2} > t_{2} - j_{3} >
t_{2}+r-j_{4}\}$ which gives 
\begin{eqnarray*}
I_{\mathcal{A},1} &=& \frac{1}{T^{2}}
\sum_{j_{1},\ldots,j_{4}=0}^{\infty}
\sum_{r,t_{1},t_{2}\in \mathcal{A}_{1}}\lambda_{M}(r)^{2}
M^{(x)}_{j_{1}}(t_{1}-j_{1})M_{j_{2}}^{(y)}(t_{1}+r-j_{2})
M^{(x)}_{j_{3}}(t_{2}-j_{3})M^{(y)}_{j_{4}}(t_{2}+r-j_{4}). 
\end{eqnarray*}
The other $23$ cases are defined similarly (different orderings of
$t_{1}-j_{1},\ldots,t_{2}+r - j_{4}$), such that we have
$I_{\mathcal{A}} = \sum_{j=1}^{24}I_{\mathcal{A}_{j}}$. We start by bounding
$I_{\mathcal{A},1}$.  Since $t_{1} - j_{1} > t_{1}+r-j_{2} > t_{2} - j_{3} >
t+r-j_{4}$, it is easy to see that 
$M^{(x)}_{j_{1}}(t_{1}-j_{1})\sum_{r < j_{2}-j_{1}}\lambda_{M}(r)^2
M_{j_{2}}^{(y)}(t_{1}+r-j_{2})
\sum_{t_{2}<
  t_{1}-j_{1}+j_{3}}M^{(x)}_{j_{3}}(t_{2}-j_{3})M^{(y)}_{j_{4}}(t_{2}+r-j_{4})$ is a martingale over $t_{1}$, 
$M_{j_{2}}^{(y)}(t_{1}+r-j_{2})\sum_{t_{2}<
t_{1}-j_{1}+j_{3}}M^{(x)}_{j_{3}}(t_{2}-j_{3})M^{(y)}_{j_{4}}(t_{2}+r-j_{4})$
is a martingale over $r$ and
$\{M^{(x)}_{j_{3}}(t_{2}-j_{3})M^{(y)}_{j_{4}}(t_{2}+r-j_{4})\}$ is a
martingale over $t_{2}$. Thus by using Burkh\"older's inequality together with 
H\"older's inequality three times, for any $q \geq 2$ we have 
\begin{eqnarray*}
\|I_{\mathcal{A},1}\|_{q} &=& 
\frac{1}{T^{2}}\sum_{j_{1},\ldots,j_{4}=0}^{\infty}\big( 
\sum_{r,t_{1},t_{2}}\lambda_{M}(r)^{2}
\|M^{(x)}_{j_{1}}(t_{1}-j_{1})\|_{4q}^{2}
\|M_{j_{2}}^{(y)}(t_{1}+r-j_{2})\|_{4q}^{2}\\
&&\|M^{(x)}_{j_{3}}(t_{2}-j_{3})\|_{4q}^{2}
\|M^{(y)}_{j_{4}}(t_{2}+r-j_{4})\|_{4q}^{2}\big)^{1/2}.
\end{eqnarray*}
Thus by using (\ref{eq:Mjb}) we have that $\|I_{\mathcal{A},1}\|_{q} = 
O(\frac{M^{1/2}}{T^{}})$ and by the same argument we have 
$I_{\mathcal{A},j} = O(\frac{M^{1/2}}{T})$ (for $2\leq j \leq 24$). Therefore, altogether this gives $\|I_{\mathcal{A}}\|_{q} = 
O(\frac{M^{1/2}}{T})$. 
We now bound $I_{\mathcal{B}}$. We first define the random variable 
\begin{eqnarray*}
&&A_{j_{1},j_{2};i}^{(x,y)}(t_{1}-j_{1}-i) = \\
&&\Ex(M_{j_{1}}^{(x)}(t_{1}-j_{1})M_{j_{2}}^{(y)}(t_{1}-j_{1})|
\mathcal{F}_{t_{1}-j_{1}-i}) - 
\Ex(M_{j_{1}}^{(x)}(t_{1}-j_{1})M_{j_{2}}^{(y)}(t_{1}-j_{1})|
\mathcal{F}_{t_{1}-j_{1}-i}). 
\end{eqnarray*}
To bound $\|A_{j_{1},j_{2};i}^{(x,y)}(t_{1}-j_{1}-i)\|_{q}$, 
we repeatedly use Ibragimov's inequality and (\ref{eq:Mjb}) to give 
\begin{eqnarray}
 \|A_{j_{1},j_{2};i}^{(x,y)}(t_{1}-j_{1}-i)\|_{q} &\leq&
2\|\Ex(\overline{M_{j_{1}}^{(x)}(t_{1}-j_{1})M_{j_{2}}^{(y)}(t_{1}-j_{1})}|
\mathcal{F}_{t_{1}-j_{1}-i})\| \nonumber\\
&\leq&C\alpha(i)\|M_{j_{1}}^{(x)}(t_{1}-j_{1})M_{j_{2}}^{(y)}(t_{1}-j_{1})\|_{q}\leq
 C\alpha(i)\alpha(j_{1})\alpha(j_{2}).\label{eq:Ajb} 
\end{eqnarray}
This gives the representation 
\begin{eqnarray*}
\overline{M_{j_{1}}^{(x)}(t_{1}-j_{1})M_{j_{2}}^{(y)}(t_{1}-j_{1})} =
\sum_{i}A_{j_{1},j_{2};i}^{(x,y)}(t_{1}-j_{1}-i).
\end{eqnarray*}
Substituting the above representation into $I_{\mathcal{B}}$ gives
\begin{eqnarray*}
I_{\mathcal{B}}&=& \frac{1}{T^{2}}\sum_{j_{1},\ldots,j_{4},i_{1},i_{2}=0}^{\infty}
\sum_{t_{1},t_{2}\in \mathcal{B}}\lambda_{M}(j_{1}-j_{2})^{2}\big[
A_{j_{1},j_{2};i_{1}}^{(x,y)}
(t_{1}-j_{1}-i)A_{j_{3},j_{4};i_{2}}^{(x,y)}(t_{2}-j_{3}-i_{2})\\
&& - \Ex(A_{j_{1},j_{2};i_{1}}^{(x,y)}
(t_{1}-j_{1}-i)A_{j_{3},j_{4};i_{2}}^{(x,y)}(t_{2}-j_{3}-i_{2}))\big] \\
 &:=& I_{\mathcal{B},1} + I_{\mathcal{B},2} + I_{\mathcal{B},3},
\end{eqnarray*}
where
\begin{eqnarray*}
I_{\mathcal{B},1} := \frac{1}{T^{2}}\sum_{j_{1},\ldots,j_{4},i_{1},i_{2}=0}^{\infty}
\sum_{t_{1}-j_{1}-i_{1}>t_{2}-j_{3}-i_{2}}\lambda_{M}(j_{1}-j_{2})^{2}
A_{j_{1},j_{2};i_{1}}^{(x,y)}
(t_{1}-j_{1}-i)A_{j_{3},j_{4};i_{2}}^{(x,y)}(t_{2}-j_{3}-i_{2}) \\
I_{\mathcal{B},2} := \frac{1}{T^{2}}\sum_{j_{1},\ldots,j_{4},i_{1},i_{2}=0}^{\infty}
\sum_{t_{1}-j_{1}-i_{1}<t_{2}-j_{3}-i_{2}}\lambda_{M}(j_{1}-j_{2})^{2}
A_{j_{1},j_{2};i_{1}}^{(x,y)}
(t_{1}-j_{1}-i)A_{j_{3},j_{4};i_{2}}^{(x,y)}(t_{2}-j_{3}-i_{2}) \\
I_{\mathcal{B},3} := \frac{1}{T^{2}}\sum_{j_{1},\ldots,j_{4},i_{1},i_{2}=0}^{\infty}
\sum_{t_{1}-j_{1}-i_{1}=t_{2}-j_{3}-i_{2}}\lambda_{M}(j_{1}-j_{2})^{2}
\overline{A_{j_{1},j_{2};i_{1}}^{(x,y)}
(t_{1}-j_{1}-i)A_{j_{3},j_{4};i_{2}}^{(x,y)}(t_{2}-j_{3}-i_{2})}.
\end{eqnarray*}
Using similar techniques to those used to bound $\|I_{\mathcal{A},1}\|_{q}$, 
Burkh\"older's and H\"older's inequalities twice on $\|I_{\mathcal{B},1}\|_{q}$, together with 
(\ref{eq:Ajb}), we obtain the bound $\|I_{\mathcal{B},1}\|_{q} = 
O(\frac{1}{T})$. A similar argument can be used for
$\|I_{\mathcal{B},2}\|_{q} = O(\frac{1}{T})$. 
To bound $\|I_{\mathcal{B},3}\|_{q}$, we need to decompose 
\begin{eqnarray*}
A_{j_{1},j_{2};i_{1}}^{(x,y)}
(t_{1}-j_{1}-i)A_{j_{3},j_{4};i_{2}}^{(x,y)}(t_{2}-j_{3}-i_{2}) - \Ex(A_{j_{1},j_{2};i_{1}}^{(x,y)}
(t_{1}-j_{1}-i)A_{j_{3},j_{4};i_{2}}^{(x,y)}(t_{2}-j_{3}-i_{2})), 
\end{eqnarray*}
into the sum of martingale differences, using this martingale
decomposition we can use the same argument as those used above to obtain 
$\|I_{\mathcal{B},3}\|= O(\frac{1}{T^{3/2}})$. Therefore,
altogether we have $\|I_{\mathcal{B}}\|_{q} = 
O(\frac{1}{T})$. Now by using similiar arguments
and repeated decompositions into martingale differences we can show that
$\|I_{\mathcal{C}}\|_{q},\ldots,\|I_{\mathcal{F}}\|_{q} = 
O(\frac{M^{1/2}}{T^{3/2}})$. Thus we have shown that $I_{\mathcal{A}}$
is the dominating term in $\mathcal{Q}_{T} - \Ex(\mathcal{Q}_{T})$. 
Since $\mathcal{S}_{T}  = \int I_{\mathcal{A}} dF_0(x) dF_0(y)$ we have obtained
the desired result. \hfill $\Box$

\vspace{3mm}
To prove the result we use the martingale central limit theorem on 
\begin{eqnarray*}
\mathcal{S}_{T} =  \frac{1}{T^2}\int \sum_{j_{1},\ldots,j_{4}=0}^{\infty}
\sum_{t_{1},r,t_{2}\in \mathcal{A}}
\lambda_{M}(r)^{2}
M^{(x)}_{j_{1}}(t_{1}-j_{1})M_{j_{2}}^{(y)}(t_{1}+r-j_{2})
M^{(x)}_{j_{3}}(t_{2}-j_{3})M^{(y)}_{j_{4}}(t_{2}+r-j_{4})dF_{0}(x)dF_{0}(y).
\end{eqnarray*}
To do this, we use the same decompositions of $I_{\mathcal{A}}$, as
that used in the proof of Theorem \ref{thm:bounds}. We set
$\mathcal{S}_{T,i} := I_{\mathcal{A},i}$, recalling
that 
\begin{eqnarray*}
\mathcal{S}_{T,i} &=& \frac{1}{T^{2}}
\sum_{j_{1},\ldots,j_{4}=0}^{\infty}
\sum_{r,t_{1},t_{2}\in \mathcal{A},i}\lambda_{M}(r)^{2}
M^{(x)}_{j_{1}}(t_{1}-j_{1})M_{j_{2}}^{(y)}(t_{1}+r-j_{2})
M^{(x)}_{j_{3}}(t_{2}-j_{3})M^{(y)}_{j_{4}}(t_{2}+r-j_{4}), 
\end{eqnarray*} 
where $\mathcal{A}_{i}$ is an ordering of 
$\{t_{1} - j_{1}, t_{1}+r-j_{2}, t_{2} - j_{3}, t_{2}+r-j_{4}\}$.
We show that $\mathcal{S}_{T,i}$ can be written as the sum of martingale 
differences. First consider $\mathcal{S}_{T,1}$, this can be written
as  $\mathcal{S}_{T,1} = \frac{1}{T^{2}}\sum_{k=1}^{T}U_{k,1}$, 
where with a change of variables we have 
\begin{eqnarray*}
U_{k,1} =
\int \sum_{j_{1}=0}^{T-k}M_{j_{1}}(k)
\sum_{j_{2},j_{3},j_{4}}\sum_{r,t_{1}\in \tilde{\mathcal{A}}_{k,1}}
\lambda_{M}(r)^2M_{j_{2}}^{(y)}(k+j_{1}+r-j_{2})
M^{(x)}_{j_{3}}(t_{2}-j_{3})M^{(y)}_{j_{4}}(t_{2}+r-j_{4})dF_{0}(x)dF_{0}(y) 
%& 2\leq k \leq T \\ 
%\end{array}
%\right..
\end{eqnarray*}
and $\tilde{\mathcal{A}}_{k,1} = \{(r,t_{2}) \textrm{ such that } (k> k+j_{1}+r-j_{2}>t_{2}-j_{3}>t_{2}+r-j_{4})\}$. 
Using a similar argument we can decompose $\mathcal{S}_{T,i}$ as $\mathcal{S}_{T,i} =
\frac{1}{T^{2}}\sum_{k=1}^{T}U_{k,i}$ (and $U_{k,i}$ is defined similar to above).  
Therefore, altogether $\mathcal{S}_{T}$ is the sum of martingale
differences, where
$\mathcal{S}_{T}  = \frac{1}{T^{2}}\sum_{k=1}^{T}\sum_{i=1}^{24}U_{k,i}$, and $\sum_{i=1}^{24}U_{k,i} \in 
\sigma(X_{k},X_{k-1},\ldots)$ are the 
martingale differences. Therefore under the conditions in Theorem \ref{thm:bounds} we have  
\begin{eqnarray*}
\mathcal{Q}_{T}  - \Ex(\mathcal{Q}_{T})= \mathcal{S}_{T}  + O_{p}(\frac{1}{T} + 
\frac{M^{1/2}}{T^{2}}) = 
\frac{1}{T^{2}}\sum_{k=1}^{T}\sum_{i=1}^{24}U_{k,i} + 
O_{p}(\frac{1}{T} + \frac{M^{1/2}}{T^{2}}). 
\end{eqnarray*}
These approximations will allow us to use the martingale 
central limit theorem to prove
asymptotic normality, which requires the following lemma.  

\begin{lemma}
\label{lemma:b4normal}
Suppose that Assumption \ref{assum:A} holds. Then for all 
$1\leq i \leq 24$ and $1\leq k \leq T$ we have
\begin{eqnarray}
\label{eq:b41}
\|\sum_{i=1}^{24}U_{k,i}\|_{q} = O(T^{1/2}M^{1/2})
\end{eqnarray}
\begin{eqnarray}
\label{eq:b42}
\frac{1}{T^{2}M}
\sum_{k=1}^{T}\Ex\big(\sum_{i=1}^{24}U_{k,i}^{2}\big) 
\rightarrow 
\frac{4}{M}\int \int
\Delta_{M}(\theta_{1} - \theta_{2})^2\prod_{i=1}^{2}
G(x_{1},x_{2};\theta_{i})G(y_{1},y_{2};\theta_{i})d\theta_{i}
\prod_{j=1}^{2}dF_{0}(x_{j})dF_{0}(y_{j})\quad
\end{eqnarray}
and 
\begin{eqnarray}
\label{eq:b43}
\frac{1}{T^{2}M}\sum_{k=1}^{T}
\bigg[\Ex\bigg((\sum_{i=1}^{24}U_{k,i})^{2}|\mathcal{F}_{k-1}\bigg)
 - \Ex\big(\sum_{i=1}^{24}U_{k,i}\big)^{2}\bigg]\Pcon 0.
\end{eqnarray}
\end{lemma}
PROOF. To prove the result we concentrate on $U_{k,1}$, a similar proof 
applies to the other terms. By using the H\"older inequality, for any
$q \geq 2$, we obtain 
\begin{eqnarray*}
\|U_{k,1}\|_{q} &\leq& \int \sum_{j_{1}=0}^{T-k}\|M_{j_{1}}(k)\|_{4q} 
\big\|\sum_{j_{2},j_{3},j_{4}}\sum_{r,t_{1}\in \mathcal{\tilde{A}}_{k,1}}
\lambda_{M}(r)^2M_{j_{2}}^{(y)}(k+j_{1}+r-j_{2})\times \\
&&
M^{(x)}_{j_{3}}(t_{2}-j_{3})M^{(y)}_{j_{4}}(t_{2}+r-j_{4})\|_{4q/3}dF_{0}(x)dF_{0}(y).
\end{eqnarray*}
Now by repeated use of Burkh\"older's inequality we have  
$\|U_{k,1}\|_{q} = O(M^{1/2}T^{1/2})$, using a similar method we obtain 
a similar bound for $\|U_{k,i}\|_{q}$, this gives (\ref{eq:b41}).

The proof of (\ref{eq:b42}) follows from the proof of Theorem \ref{thm:bounds}
(noting that the asymptotic variance of $\mathcal{Q_{T}}$ is determined by the 
variance of $\mathcal{S}_{T}$).

To prove (\ref{eq:b43}), we consider only the $U_{k,1}$ (the proof
involving the other terms in similar). For brevity we
write $U_{k,1}$ as
\begin{eqnarray*}
U_{k,1} = \int \sum_{j_{1}=0}^{T-k}M_{j_{1}}^{(x)}(k)
N_{j_{1},k-1,1}^{(x,y)}dF_{0}(x)dF_{0}(y), 
\end{eqnarray*}
where 
\begin{eqnarray*}
N_{j_{1},k-1,1}^{(x,y)} &=& \sum_{j_{2},j_{3},j_{4}}
\sum_{r,t_{1}\in \tilde{\mathcal{A}}_{k,1}}
\lambda_{M}(r)^2M_{j_{2}}^{(y)}(k+j_{1}+r-j_{2})
M^{(x)}_{j_{3}}(t_{2}-j_{3})M^{(y)}_{j_{4}}(t_{2}+r-j_{4}).
\end{eqnarray*}
Noting that $N_{j_{1},k-1,1}^{(x,y)}\in\mathcal{F}_{k-1}$ we have 
\begin{eqnarray*}
&&\hskip -1em \frac{1}{T^{2}M}\sum_{k=1}^{T}\big(\Ex\big(U_{k,1})^{2}|\mathcal{F}_{k-1}\big)
 - \Ex\big(U_{k,1}\big)^{2}\big) = \\
&&\hskip -1em \frac{1}{T^{2}M}\sum_{k=1}^{T}
\int \sum_{j_{1},j_{2}=0}^{T-k}
\big(\Ex(M_{j_{1}}^{(x)}(k)M_{j_{2}}^{(x_{1})}(k)|\mathcal{F}_{k-1}) -
\Ex(M_{j_{1}}^{(x)}(k)M_{j_{2}}^{(x_{2})}(k))\big)
N_{j_{1},k-1,1}^{(x_{1},y_{2})}N_{j_{1},k-1,1}^{(x_{2},y_{2})}
\prod_{i=1}^{2}dF_{0}(x_{i})dF_{0}(y_{i}) \\
&\hskip -1em  +& \hskip -1em \frac{1}{T^{2}M}\sum_{k=1}^{T}
\int \sum_{j_{1},j_{2}=0}^{T-k}
\Ex(M_{j_{1}}^{(x)}(k)M_{j_{2}}^{(x_{2})}(k))\big)
\big(
N_{j_{1},k-1,1}^{(x_{1},y_{2})}N_{j_{1},k-1,1}^{(x_{2},y_{2})} - 
\Ex(N_{j_{1},k-1,1}^{(x_{1},y_{2})}N_{j_{1},k-1,1}^{(x_{2},y_{2})})\big)
\prod_{i=1}^{2}dF_{0}(x_{i})dF_{0}(y_{i}). 
\end{eqnarray*}
Now by using similar methods to the iterative martingale methods detailed in
the proof of Theorem \ref{thm:bounds}, we can show that the
$\|\cdot\|_{q}$-norm ($q\geq 2$) of the above converges to zero, thus we have (\ref{eq:b43}).

\vspace{3mm}

{\bf PROOF of Theorem \ref{theorem:null}}  Using the above we have 
$\mathcal{Q}_{T} =
\frac{1}{T^{2}}\sum_{k=1}^{T}\sum_{i=1}^{24}U_{k,i} + 
O_{p}(\frac{1}{T} + \frac{M^{1/2}}{T^{2}})$, thus $\mathcal{Q}_{T} -
\Ex(\mathcal{Q}_{T})$ can be written
as the sum of martingales plus a smaller order term. 
Therefore to prove asymptotic normality
of $\mathcal{Q}_{T}$ we can use the martingale central limit, for this  
we need to verify (a) the conditional Lindeberg condition 
$\frac{1}{MT^{2}}\sum_{k=1}^{T}
\Ex(|\sum_{i=1}^{24}U_{k,i}|^{2}I(\frac{1}{M^{1/2}T}
|\sum_{i=1}^{24}U_{k,i}| > \varepsilon)|\mathcal{F}_{k-1})\Pcon 0$ 
for all $\varepsilon > 0$, (b) 
that $\frac{1}{MT^{2}}\sum_{k=1}^{T}\Ex(|\sum_{i=1}^{24}U_{k,i}|^{2}|
\mathcal{F}_{k-1})-\frac{T^{2}}{M}\var(Q_{T})\Pcon 0$.  

To verify the conditional Lindeberg condition, we observe that the Cauchy-Schwartz and
Markov's inequalities give 
\begin{eqnarray*}
\frac{1}{MT^{2}}\sum_{k=1}^{T}
\Ex\big(|\sum_{i=1}^{24}U_{k,i}|^{2}I(\frac{1}{M^{1/2}T}
|\sum_{i=1}^{24}U_{k,i}| > \varepsilon)|\mathcal{F}_{k-1}\big) \leq 
\frac{1}{\varepsilon M^{2}T^{4}}\sum_{k=1}^{T}
\Ex(|\sum_{i=1}^{24}U_{k,i}|^{4}|\mathcal{F}_{k-1}):= B_{T}.
\end{eqnarray*}
By using (\ref{eq:b41}) the expectation of the above is 
$\Ex(B_{T}) = O(\frac{1}{T})$. As $B_{T}$ is a  non-negative random variable, this implies 
$B_{T}\Pcon 0$ as $T\rightarrow \infty$. Thus we have shown that the Lindeberg 
condition is satisfied. To prove (b) we note that 
\begin{eqnarray*}
&&\frac{1}{MT^{2}}\sum_{k=1}^{T}\Ex\bigg(|\sum_{i=1}^{24}U_{k,i}|^{2}|
\mathcal{F}_{k-1}\bigg)  - \frac{T^{2}}{M}\var(Q_{T}) \\
&=&\frac{1}{MT^{2}}\sum_{k=1}^{T}\bigg[\Ex\big(|\sum_{i=1}^{24}U_{k,i}|^{2}|
\mathcal{F}_{k-1}\big) - \Ex(|\sum_{i=1}^{24}U_{k,i}|^{2})\bigg]  + 
 \frac{1}{MT^{2}}\sum_{k=1}^{T}
\Ex(|\sum_{i=1}^{24}U_{k,i}|^{2}) - \frac{T^{2}}{M}\var(Q_{T}).
\end{eqnarray*}
By using (\ref{eq:b42}) and (\ref{eq:b43}) the above converges to zero in 
probability. Thus we have verified the conditions of the martingale central
limit theorem and we have the desired result. \hfill $\Box$

\vspace{3mm}

\subsection{Proof of Theorem \ref{theorem:alternative}}

As the limiting distribution of $\mathcal{Q}_{T}$ is determined by
$\mathcal{Q}_{T,2}$, we rewrite $\mathcal{Q}_{T,2}$ in such a way that
the same methods used to prove Theorem 2 in \citeA{p:lee-sub-11}, can be
used to obtain the limiting distribution. We observe that
\begin{eqnarray*}
\mathcal{Q}_{T,2} &=& \frac{2}{T}\Re\int \sum_{k} \Lambda_{T}(x,y;\omega_{k})
\big\{J_{T}(x;\omega_{k})\overline{J}_{T}(y;\omega_{k}) -
\Ex(J_{T}(x;\omega_{k})\overline{J}_{T}(y;\omega_{k}))\big\}dF_{0}(x)dF_{0}(y)\\ 
&=& \int \frac{2}{T}\sum_{t,\tau}\lambda_{M}(t-\tau)^{2}D_{t-\tau,T}(x,y)
(Z_{t}(x)Z_{\tau}(y) - \Ex(Z_{t}(x)Z_{\tau}(y)))dF_{0}(x)dF_{0}(y) \nonumber\\
&=&  \int \frac{2}{T}\sum_{t,\tau}\lambda_{M}(t-\tau)^{2}D_{t-\tau,T}(x,y)
(\tilde{Z}_{T}(x)\tilde{Z}_{\tau}(y) - C_{r}(x,y))dF_{0}(x)dF_{0}(y) +  
O_{p}(\frac{1}{T}), \label{eq:simple}
\end{eqnarray*}
where $\Lambda_{T}(x,y;\omega_{s}) = 
\sum_{r}\lambda_{M}(r)^{2}(\frac{T-|r|}{T})[C_{r,1}(x,y) -C_{r,0}(x,y)]\exp(ir\omega_{k})$, 
$D_{r,T}(x,y) = (\frac{T-|r|}{T})
[C_{r,1}(x,y) - C_{r,0}(x,y)]$ and 
$\tilde{Z}_{t}(x) = I(X_{t} \leq x) - F_{1}(x)$. 

\vspace{3mm}

{\bf PROOF of Theorem \ref{theorem:alternative}} 
Now we observe that under the stated assumptions of the theorem we
have  that the quantile covariances under the null decay at the rate
$\sup_{x,y}|C_{r,0}(x,y)|\leq K|r|^{-(2+\delta)}$ (for some $\delta > 0$) and 
$\sup_{x,y}|C_{r,1}(x,y)|\leq K|r|^{-s}$ (for some $s
> 2$). Thus by definition of $D_{r,T}(\cdot)$, we have 
$\sup_{x,y}|\lambda_{M}(r)D_{r,T}(x,y)|\leq K|r|^{-\min(2+\delta,s)}$.  
Thus we can write $\mathcal{Q}_{T,2}$ as
\begin{eqnarray*}
\mathcal{Q}_{T,2} 
&=&  \int
\frac{2}{T}\sum_{t,\tau}\lambda_{M}(t-\tau)^{2}D_{t-\tau,T}(x,y)
(\tilde{Z}_{T}(x)\tilde{Z}_{\tau}(y) - \Ex(\tilde{Z}_{T}(x)\tilde{Z}_{\tau}(y)))dF_{0}(x)dF_{0}(y) +  
O_{p}(\frac{1}{T}),
\end{eqnarray*}
where we observe that terms where $|t - \tau| > 2M$, are zero. 
Thus using the Bernstein blocking arguments for quadratic forms used to prove Theorem 
2, Lee and Subba Rao (2011), \nocite{p:lee-sub-11} we can show
asymptotic normality of the above.
This proves (\ref{eq:Q2}).  Finally to prove (\ref{eq:Q3}), we note
that $\mathcal{Q}_{T} = \mathcal{Q}_{T,2} + E_{T,2} +
O_{p}(\frac{M^{1/2}}{T} + \frac{M}{T} + \frac{1}{M^{s-1}})$, by using
(\ref{eq:Q2}), this immediately leads to (\ref{eq:Q3}).  
\hfill $\Box$

%%% Local Variables: 
%%% mode: latex
%%% TeX-master: t
%%% End: 

\bibliographystyle{apacite}
\bibliography{goodness_of_fit}

\begin{table}[ht]
\caption{$H_0: AR(1)$, $H_{A}:ARCH$ $T=100$}
\begin{center}
\label{tab:1}
\begin{tabular}{c!{\VRule[0.4pt]}r!{\VRule[0.8pt]}r!{\VRule[0.4pt]}r!{\VRule[0.4pt]}r!{\VRule[0.4pt]}r!{\VRule[0.8pt]}r!{\VRule[0.4pt]}r!{\VRule[0.4pt]}r!{\VRule[0.4pt]}r}
\specialrule{0.8pt}{0pt}{0pt}
\multicolumn{2}{c!{\VRule[0.8pt]}}{\multirow{2}{*}{ $T=100$} }
& \multicolumn{4}{c!{\VRule[0.8pt]}}% 
{\centering $\alpha = 0.1$  }& \multicolumn{4}{c}% 
{\centering $\alpha = 0.05$  }\\
\cline{3-10}
\multicolumn{2}{c!{\VRule[0.8pt]}}{} & \multicolumn{2}{c!{\VRule[0.4pt]}}{Bootstrap} & \multicolumn{2}{c!{\VRule[0.8pt]}}{Normal}& \multicolumn{2}{c!{\VRule[0.4pt]}}{Bootstrap} & \multicolumn{2}{c}{Normal}\\
\hline
\multicolumn{1}{c!{\VRule[0.4pt]}}{$a$} & \multicolumn{1}{c!{\VRule[0.8pt]}}{M} & \multicolumn{1}{c!{\VRule[0.4pt]}}{$H_0$} & \multicolumn{1}{c!{\VRule[0.4pt]}}{$H_A$} & \multicolumn{1}{c!{\VRule[0.4pt]}}{$H_0$} & \multicolumn{1}{c!{\VRule[0.8pt]}}{$H_A$} & \multicolumn{1}{c!{\VRule[0.4pt]}}{$H_0$} & \multicolumn{1}{c!{\VRule[0.4pt]}}{$H_A$} & \multicolumn{1}{c!{\VRule[0.4pt]}}{$H_0$} & \multicolumn{1}{c}{$H_A$} \\
\specialrule{0.8pt}{0pt}{0pt}
\multirow{4}{*}{0.3}&11&0.052&1&0.076&1&0.021&0.972&0.054&1\\
&16&0.04&0.869&0.062&0.971&0.011&0.262&0.04&0.854\\
&21&0.048&0.386&0.064&0.561&0.021&0.106&0.043&0.348\\
&25&0.021&0.071&0.048&0.229&0.014&0.016&0.029&0.12\\
\hline\multirow{4}{*}{0.4}&11&0.048&1&0.082&1&0.02&1&0.055&1\\
&16&0.043&1&0.059&1&0.013&0.939&0.041&1\\
&21&0.046&0.932&0.066&0.997&0.011&0.416&0.046&0.929\\
&25&0.036&0.582&0.055&0.832&0.01&0.124&0.037&0.598\\
\hline\multirow{4}{*}{0.5}&11&0.046&1&0.073&1&0.015&1&0.052&1\\
&16&0.049&1&0.078&1&0.027&1&0.045&1\\
&21&0.046&1&0.06&1&0.015&0.985&0.037&1\\
&25&0.047&1&0.062&1&0.015&0.397&0.043&1\\
\hline\multirow{4}{*}{0.55}&11&0.041&1&0.096&1&0.018&1&0.057&1\\
&16&0.045&1&0.066&1&0.017&1&0.046&1\\
&21&0.065&1&0.06&1&0.034&1&0.034&1\\
&25&0.045&1&0.051&1&0.024&1&0.032&1\\
\specialrule{0.8pt}{0pt}{0pt}
\end{tabular} 
\end{center}
\end{table}

\begin{table}
\caption{$H_0: AR(1)$, $H_{A}:ARCH(1)$ $T = 500$}
\begin{center}
\label{tab:2}
\begin{tabular}{c!{\VRule[0.4pt]}r!{\VRule[0.8pt]}r!{\VRule[0.4pt]}r!{\VRule[0.4pt]}r!{\VRule[0.4pt]}r!{\VRule[0.8pt]}r!{\VRule[0.4pt]}r!{\VRule[0.4pt]}r!{\VRule[0.4pt]}r}
\specialrule{0.8pt}{0pt}{0pt}
\multicolumn{2}{c!{\VRule[0.8pt]}}{\multirow{2}{*}{ $T=500$} }
& \multicolumn{4}{c!{\VRule[0.8pt]}}% 
{\centering $\alpha = 0.1$  }& \multicolumn{4}{c}% 
{\centering $\alpha = 0.05$  }\\
\cline{3-10}
\multicolumn{2}{c!{\VRule[0.8pt]}}{} & \multicolumn{2}{c!{\VRule[0.4pt]}}{Bootstrap} & \multicolumn{2}{c!{\VRule[0.8pt]}}{Normal}& \multicolumn{2}{c!{\VRule[0.4pt]}}{Bootstrap} & \multicolumn{2}{c}{Normal}\\
\hline
\multicolumn{1}{c!{\VRule[0.4pt]}}{$a$} & \multicolumn{1}{c!{\VRule[0.8pt]}}{M} & \multicolumn{1}{c!{\VRule[0.4pt]}}{$H_0$} & \multicolumn{1}{c!{\VRule[0.4pt]}}{$H_A$} & \multicolumn{1}{c!{\VRule[0.4pt]}}{$H_0$} & \multicolumn{1}{c!{\VRule[0.8pt]}}{$H_A$} & \multicolumn{1}{c!{\VRule[0.4pt]}}{$H_0$} & \multicolumn{1}{c!{\VRule[0.4pt]}}{$H_A$} & \multicolumn{1}{c!{\VRule[0.4pt]}}{$H_0$} & \multicolumn{1}{c}{$H_A$} \\
\specialrule{0.8pt}{0pt}{0pt}
\multirow{4}{*}{0.3}&14&0.053&1&0.098&1&0.024&1&0.063&1\\
&21&0.064&1&0.082&1&0.023&1&0.052&1\\
&28&0.06&1&0.093&1&0.024&1&0.062&1\\
&35&0.07&1&0.086&1&0.033&1&0.062&1\\
\hline\multirow{4}{*}{0.4}&14&0.043&1&0.092&1&0.014&1&0.064&1\\
&21&0.058&1&0.092&1&0.015&1&0.056&1\\
&28&0.066&1&0.094&1&0.03&1&0.061&1\\
&35&0.073&1&0.087&1&0.032&1&0.052&1\\
\hline\multirow{4}{*}{0.5}&14&0.031&1&0.105&1&0.018&1&0.072&1\\
&21&0.059&1&0.079&1&0.03&1&0.05&1\\
&28&0.076&1&0.111&1&0.046&1&0.069&1\\
&35&0.053&1&0.086&1&0.022&1&0.055&1\\
\hline\multirow{4}{*}{0.55}&14&0.038&1&0.107&1&0.014&1&0.077&1\\
&21&0.056&1&0.108&1&0.021&1&0.067&1\\
&28&0.071&1&0.103&1&0.032&1&0.06&1\\
&35&0.051&1&0.089&1&0.026&1&0.06&1\\
\specialrule{0.8pt}{0pt}{0pt}
\end{tabular} 
\end{center}
\end{table}

\begin{table}
\caption{$H_0: ARCH(1)$, $H_{A} = AR(1)$ $T=100$}
\begin{center}
\label{tab:3}
\begin{tabular}{c!{\VRule[0.4pt]}r!{\VRule[0.8pt]}r!{\VRule[0.4pt]}r!{\VRule[0.4pt]}r!{\VRule[0.4pt]}r!{\VRule[0.8pt]}r!{\VRule[0.4pt]}r!{\VRule[0.4pt]}r!{\VRule[0.4pt]}r}
\specialrule{0.8pt}{0pt}{0pt}
\multicolumn{2}{c!{\VRule[0.8pt]}}{\multirow{2}{*}{ $T=100$} }
& \multicolumn{4}{c!{\VRule[0.8pt]}}% 
{\centering $\alpha = 0.1$  }& \multicolumn{4}{c}% 
{\centering $\alpha = 0.05$  }\\
\cline{3-10}
\multicolumn{2}{c!{\VRule[0.8pt]}}{} & \multicolumn{2}{c!{\VRule[0.4pt]}}{Bootstrap} & \multicolumn{2}{c!{\VRule[0.8pt]}}{Normal}& \multicolumn{2}{c!{\VRule[0.4pt]}}{Bootstrap} & \multicolumn{2}{c}{Normal}\\
\hline
\multicolumn{1}{c!{\VRule[0.4pt]}}{$a$} & \multicolumn{1}{c!{\VRule[0.8pt]}}{M} & \multicolumn{1}{c!{\VRule[0.4pt]}}{$H_0$} & \multicolumn{1}{c!{\VRule[0.4pt]}}{$H_A$} & \multicolumn{1}{c!{\VRule[0.4pt]}}{$H_0$} & \multicolumn{1}{c!{\VRule[0.8pt]}}{$H_A$} & \multicolumn{1}{c!{\VRule[0.4pt]}}{$H_0$} & \multicolumn{1}{c!{\VRule[0.4pt]}}{$H_A$} & \multicolumn{1}{c!{\VRule[0.4pt]}}{$H_0$} & \multicolumn{1}{c}{$H_A$} \\
\specialrule{0.8pt}{0pt}{0pt}
\multirow{4}{*}{0.3}&11&0.039&0.994&0.08&0.997&0.022&0.984&0.051&0.995\\
&16&0.043&0.978&0.086&0.991&0.009&0.925&0.055&0.983\\
&21&0.045&0.98&0.07&0.99&0.016&0.934&0.051&0.983\\
&25&0.026&0.939&0.059&0.976&0.011&0.895&0.045&0.965\\
\hline\multirow{4}{*}{0.4}&11&0.046&1&0.086&1&0.012&0.999&0.053&1\\
&16&0.049&0.993&0.092&0.999&0.014&0.988&0.062&0.996\\
&21&0.03&0.994&0.07&0.997&0.017&0.983&0.046&0.997\\
&25&0.038&0.994&0.083&0.997&0.024&0.982&0.059&0.994\\
\hline\multirow{4}{*}{0.5}&11&0.054&1&0.107&1&0.024&1&0.067&1\\
&16&0.063&1&0.098&1&0.03&1&0.066&1\\
&21&0.051&1&0.083&1&0.022&1&0.061&1\\
&25&0.028&0.997&0.06&0.998&0.012&0.995&0.043&0.998\\
\hline\multirow{4}{*}{0.55}&11&0.074&1&0.113&1&0.03&1&0.081&1\\
&16&0.056&1&0.087&1&0.02&1&0.054&1\\
&21&0.065&1&0.08&1&0.038&1&0.057&1\\
&25&0.067&1&0.088&1&0.03&1&0.065&1\\
\specialrule{0.8pt}{0pt}{0pt}
\end{tabular} 
\end{center}
\end{table}

\begin{table}
\caption{$H_0: ARCH(1)$ and $H_{A}:AR(1)$ $T=500$}
\begin{center}
\label{tab:4}
\begin{tabular}{c!{\VRule[0.4pt]}r!{\VRule[0.8pt]}r!{\VRule[0.4pt]}r!{\VRule[0.4pt]}r!{\VRule[0.4pt]}r!{\VRule[0.8pt]}r!{\VRule[0.4pt]}r!{\VRule[0.4pt]}r!{\VRule[0.4pt]}r}
\specialrule{0.8pt}{0pt}{0pt}
\multicolumn{2}{c!{\VRule[0.8pt]}}{\multirow{2}{*}{ $T=500$} }
& \multicolumn{4}{c!{\VRule[0.8pt]}}% 
{\centering $\alpha = 0.1$  }& \multicolumn{4}{c}% 
{\centering $\alpha = 0.05$  }\\
\cline{3-10}
\multicolumn{2}{c!{\VRule[0.8pt]}}{} & \multicolumn{2}{c!{\VRule[0.4pt]}}{Bootstrap} & \multicolumn{2}{c!{\VRule[0.8pt]}}{Normal}& \multicolumn{2}{c!{\VRule[0.4pt]}}{Bootstrap} & \multicolumn{2}{c}{Normal}\\
\hline
\multicolumn{1}{c!{\VRule[0.4pt]}}{$a$} & \multicolumn{1}{c!{\VRule[0.8pt]}}{M} & \multicolumn{1}{c!{\VRule[0.4pt]}}{$H_0$} & \multicolumn{1}{c!{\VRule[0.4pt]}}{$H_A$} & \multicolumn{1}{c!{\VRule[0.4pt]}}{$H_0$} & \multicolumn{1}{c!{\VRule[0.8pt]}}{$H_A$} & \multicolumn{1}{c!{\VRule[0.4pt]}}{$H_0$} & \multicolumn{1}{c!{\VRule[0.4pt]}}{$H_A$} & \multicolumn{1}{c!{\VRule[0.4pt]}}{$H_0$} & \multicolumn{1}{c}{$H_A$} \\
\specialrule{0.8pt}{0pt}{0pt}
\multirow{4}{*}{0.3}&14&0.072&1&0.09&1&0.025&1&0.059&1\\
&21&0.062&1&0.094&1&0.032&1&0.059&1\\
&28&0.067&1&0.097&1&0.024&1&0.062&1\\
&35&0.076&1&0.101&1&0.026&1&0.073&1\\
\hline\multirow{4}{*}{0.4}&14&0.045&1&0.097&1&0.022&1&0.059&1\\
&21&0.075&1&0.105&1&0.03&1&0.077&1\\
&28&0.06&1&0.111&1&0.024&1&0.07&1\\
&35&0.085&1&0.12&1&0.041&1&0.086&1\\
\hline\multirow{4}{*}{0.5}&14&0.053&1&0.129&1&0.032&1&0.079&1\\
&21&0.1&1&0.121&1&0.054&1&0.082&1\\
&28&0.111&1&0.124&1&0.071&1&0.085&1\\
&35&0.066&1&0.117&1&0.029&1&0.075&1\\
\hline\multirow{4}{*}{0.55}&14&0.099&1&0.143&1&0.047&1&0.104&1\\
&21&0.074&1&0.119&1&0.042&1&0.083&1\\
&28&0.078&1&0.11&1&0.037&1&0.072&1\\
&35&0.082&1&0.119&1&0.037&1&0.085&1\\
\specialrule{0.8pt}{0pt}{0pt}
\end{tabular} 
\end{center}
\end{table}

%%% Local Variables: 
%%% mode: latex
%%% TeX-master: t
%%% End: 

\end{document}